\input amstex
\documentstyle{amsppt}
\def\sideremark#1{\relax}
\define\ups{{\text{\rm \char '007}}}
\let\ccdot\cdot
\def\cdot{\hbox to 2.5pt{\hss$\ccdot$\hss}}
\define\ox{\otimes}
\define\gog{{\goth g}}
\define\tsum{\tsize\sum}
\define\newquad{\hbox{\hphantom{${}={}$}}}
\define\na{\nabla}
\define\>{\rightarrow}
\define\<{\leftarrow}
\define\[{\lbrack}
\define\]{\rbrack}
\redefine\o{\circ}

\define\al{\alpha}
\define\be{\beta}
\define\ga{\gamma}
\define\de{\delta}

\define\ze{\zeta}
\define\et{\eta}
\define\th{\theta}

\define\ka{\kappa}
\define\la{\lambda}
\define\rh{\rho}
\define\si{\sigma}
\define\ta{\tau}
\define\ph{\varphi}
\define\ch{\chi}
\define\ps{\psi}
\define\om{\omega}
\define\Ga{\Gamma}

\define\La{\Lambda}

\define\Ph{\Phi}

\define\Om{\Omega}

\define\ad{{\operatorname{ad}}}
\define\Ad{{\operatorname{Ad}}}
\let\ii\i
\redefine\i{^{-1}}
\define\row#1#2#3{#1_{#2},\ldots,#1_{#3}}
\define\rowup#1#2#3{#1^{#2},\ldots,#1^{#3}}
\define\irow#1#2#3{{#1_{#2}\ldots#1_{#3}}}
\define\irowup#1#2#3{{#1^{#2}\ldots#1^{#3}}}
\define\x{\times}
\define\nmb#1#2{#2}      
\def\idx{}               

\define\C{{\Cal C}}
\define\Mf{\Cal Mf}

\define\J{{\Cal J}}


\redefine\L{{\Cal L}}
\predefine\SS\S
\redefine\S{{\Cal S}}

\define\dd#1{\tfrac \partial{\partial #1}}
\define\field#1#2#3{{#1}^{#2}\frac\partial{\partial{#1}_{#3}}}

\define\pol#1 #2#3{L^{#1}_{\text{sym}}(\Bbb R^#2,\Bbb R^#3)}
\define\sempr{\rtimes}

\def\today{\ifcase\month\or
 January\or February\or March\or April\or May\or June\or
 July\or August\or September\or October\or November\or December\fi
 \space\number\day, \number\year}
\topmatter
\title Invariant operators on manifolds with almost Hermitian
symmetric structures, \\I. Invariant differentiation\endtitle
\author Andreas \v Cap, Jan Slov\'ak, Vladim\'ir Sou\v cek\endauthor
\abstract
This is the first part of a series of papers. The whole series aims to
develop the tools for the study of all almost Hermitian symmetric structures
in a unified way. In particular, methods for the construction of invariant
operators, their classification and the study of their properties
will be worked out.

In this paper we present the invariant
differentiation with respect to a Cartan connection and we expand
the differentials in the terms of the underlying linear
connections belonging to the structures in question. Then we discuss the
holonomic and non-holonomic jet extensions and we suggest methods for
the construction of invariant operators.
\endabstract
\address
Institut f\"ur Mathematik, Universit\"at Wien, Strudlhofgasse 4,
1090 Wien, Austria\newline\indent
Department of Algebra and Geometry, Masaryk University in Brno,
Jan\'a\v ckovo n\' am\. 2a, 662~95 Brno, Czech
Republic\newline\indent
Mathematical Institute, Charles University, Sokolovsk\'a 83,
Praha, Czech Republic
\endaddress
\thanks This work was mostly done during the stay of the authors at the
Erwin Schr\"odinger International Institute of Mathematical Physics in
Vienna. The second and the third authors are also supported by the
GA\v CR, grant Nr\. 2178
\endthanks
\rightheadtext{I. Invariant differentiation}
\endtopmatter

\document
\head\nmb0{1}. Introduction \endhead

It is well known that the theories of conformal Riemannian structures and
of projective structures admits a unified exposition in terms of the so
called $|1|$-graded Lie algebras $\goth g=\goth g_{-1}\oplus \goth g_0\oplus
\goth g_1$, see e.g\. \cite{Kobayashi, 72}, and there has been a wide
discussion on geometries fitting into a similar scheme, see e.g\.
\cite{Kobayashi, Nagano, 64, 65}, \cite{Ochiai, 70}. Already there, the
Cartan connections appeared as the absolute parallelisms obtained on the
last non-trivial prolongations of the original $G$-structure in question, and
it turned out that they should play a role similar to that of the Riemannian
connections in Riemannian geometry. This was the point of view adopted by
\cite{Ochiai} under the strong additional conditions of the vanishing of the
torsions. A different approach covering all $|1|$-graded algebras
$\goth g$ can be found in \cite{Baston, 91}.

In our setting, we shall consider a fixed connected Lie group $G$ with such an
$|1|$-graded Lie algebra $\goth g$, its subgroup $B$ with the Lie algebra
$\goth b=\goth g_0\oplus \goth g_1$, the normal subgroup $B_1\subset B$ with
the Lie algebra $\goth g_1$ and the Lie group $B_0=B/B_1$ with the Lie
algebra $\goth g_0$. The corresponding geometric structures are then
reductions of the linear frame bundles $P^1M$ on $\operatorname{dim}\goth
g_{-1}$--dimensional manifolds $M$ to the structure group $B_0$. It turns
out that the flat (homogeneous) models for such structures are the Hermitian
symmetric spaces $G/B$ and, following \cite{Baston, 91}, we call them the
almost Hermitian symmetric structures, briefly the AHS structures.
Similar structures were studied earlier by \cite{Goncharov,
1987}.

Our present goal is to discuss the role played by the Cartan connections in
the general theory. More explicitly, we aim to develop a calculus for the
Cartan connections similar to the Ricci calculus for the linear (Riemannian)
connections. Thus we postpone the general construction of the canonical
Cartan connections to the next part in the series, while now we discuss the
AHS structures in a more abstract form, as principal $B$-bundles equipped
with an analogy to the soldering form on the linear frame bundles. This
corresponds to thinking about conformal and projective structures as being
second order structures (i.e\. reductions of the second order frame bundles
$P^2M$). In particular, we work out the tools for building the invariant
operators as expressions in terms of the linear connections belonging to the
structures and we describe an alternative to the jet extensions, the
semi-holonomic jet extensions. The  $k$th order jet extensions fail to be
associated bundles to the defining principal $B$-bundles, except for the
locally flat structures, however the semi-holonomic ones always
are.  We construct
a universal invariant differential operator with values in the
semi-holonomic jets, the
iterated invariant derivative with respect to a Cartan connection. This
approach generalizes vastly most of the classical constructions of invariant
operators in the conformal Riemannian geometries. We shall comment more
explicitly on the direct links in the text.

The next part of the series shows that all AHS structures,
defined as first order structures in the way indicated above, give rise to
canonical principal $B$-bundles equipped with the canonical soldering forms.
Moreover we shall give an explicit construction of the canonical Cartan
connections there. Thus, the calculus developed here suggests direct
methods for the study of invariant operators on all AHS structures.

In the third part, we shall rewrite the recurrence procedure for the
expansion of the invariant differentials in the terms of the finite
dimensional representation theory of the invidual Lie algebras. This will
help us to achieve an explicit construction of large classes of invariant
operators, even for the `curved cases'. In particular, we shall show that
all the operators on the locally flat AHS manifolds, well known from the
theory of the generalized Verma modules as the standard operators, have a
canonical extension to all AHS manifolds. Moreover, we shall even present
universal formulae for those operators in a closed form in terms of the
linear connections belonging to the structure.

The inspiration for the fourth part comes from the Jantzen-Zuckerman
translation functors on the generalized Verma modules which present the key
to the classification of all the operators on the locally flat AHS
manifolds. Several attempts to find an analogy of this translation principle
for the conformal Riemannian geometries, working also outside of the localy
flat spaces, can be found in the literature. The first one was
\cite{Eastwood, Rice, 87}, see the survey on conformally invariant operators
\cite{Baston, Eastwood, 90} for more information. The only attempt to work
out the translation procedure for all AHS structures was done in the second
part of \cite{Baston, 91}, however the argumentation is very opaque
there. In
our approach, we shall obtain a complete and explicit version of the
translation principle by analyzing the algebraic structure of the
semi-holonomic jets, in a very straightforward way. In fact our arguments
will follow the lines of the Eastood's `curved translation principle' in the
conformal Riemannian geometry, see e.g\. \cite{Eastwood, 95}.

More detailed links to the existing literature will be given at the suitable
places in the individual parts of this series. But let us say at least, that
our motivation comes mostly from the wide range of results on the invariant
operators on conformal Riemannian manifolds, in particular the
series of papers by T.~N.~Bailey, R.~J.~Baston, T.~P.~Branson,
M.~G.~Eastwood, C.~R.~Graham, H.~P.~Jakobsen, V.~W\"unsch and others, cf\.
the references at the end of this paper. Our development is probably most
influenced by \cite{Baston, 90, 91}.

Let us conclude the Introduction with a brief description of the structure of
this paper. First we study the obvious operation on the frame forms of
sections of associated bundles defined by means of the horizontal vector
fields with respect to the Cartan connections. However we exploit the very
special properties of the structures and connections in question, and we can
iterate our derivatives.  The result of such an iterated differentiation of
a section is not the frame form of a section in general (i.e\. the required
equivariance properties fail).

In the next two sections, we develop a
recurrence procedure which enables us to get expressions for the iterated
derivatives in terms of the linear connections belonging to the structure
and to localize the failure of their equivariance.  Once we will have the
canonical Cartan connections, this calculus can be used directly for the
classification in low orders but also for direct construtions of higher
order operators, in a way generalizing vastly the approach by \cite{Fegan,
76}, \cite{Branson, 89}, \cite{W\"unsch, 86}.

In the fifth section, we work out the proper bundles, which are codomains of
the iterated differentials, the semi-holonomic jet prolongations, and we
discuss the `universal invariant differential operators' defined by means of
our invariant derivatives.  Furthermore we show, that each iterated
invariant differential splits into three parts: the usual iterated covariant
differential with respect to a linear connection belonging to the structure,
the `correction terms', and the `obstruction terms'. We also succeed in
localization of the parts of the general obstruction terms which are the
proper obstructions (i.e\. the other ones vanish as a consequence). In
particular, we find a much smaller `algebraic obstruction', the algebraic
vanishing of which suffices.

In the last section, we illustrate our methods on examples. First we
recover the canonical Cartan connection on the conformal structures by
computing explicitly the necessary deformation and then we show how our
formulae lead to some known operators.

\specialhead Acknowledgment\endspecialhead
The  work on this series of papers has been influenced by discussions with
several mathematicians, the authors like to mention especially the fruitful
communication with J.~Bure\v s and M.~Eastwood.

\head \nmb0{2}. The invariant differentiation\endhead

\subhead \nmb.{2.1}. Cartan connections\endsubhead
Let $G$ be a Lie group, $B\subset G$ a closed subgroup, and let $\goth
g$, $\goth b$ be the Lie algebras of $G$ and $B$. Further, let
$P\to M$ be a principal fiber bundle with structure group $B$ and
let us denote by $\ze_X$ the fundamental vector field
corresponding to $X\in \goth b$. A
$\goth g$-valued one form $\om\in \Om^1(P,\goth g)$ with the
properties
\roster
\item $\om(\ze_X) = X$ for all $X\in \goth b$
\item $(r^b)^*\om=\text{Ad}(b^{-1})\o \om$ for all $b\in B$
\item $\om|_{T_uP}\: T_uP\to \goth g$ is a bijection for all
$u\in P$
\endroster
is called a \idx{\it Cartan connection}. Clearly,
$\text{dim}M=\text{dim}G-\text{dim}B=\text{dim}(G/B)$ if a Cartan
connection exists.

The curvature $K\in \Om^2(P,\goth g)$ of a Cartan connection $\om$
is defined by the structure equation
$$
d\om = -\frac12[\om,\om] + K.
$$
The Cartan connection $\om$ defines for each element $Y\in \goth
g$ the vector field $\om^{-1}(Y)$ given by the equality
$\om(\om^{-1}(Y)(u))=Y$ for all $u\in P$. This defines an
absolute parallelism on $P$.

{}From now on we assume that there is an abelian subalgebra $\goth g_{-1}$
in $\goth g$ which is complementary to $\goth b$, so that
$\goth g=\goth g_{-1}\oplus \goth b$. Then
$\om^{-1}(\goth g_{-1})\subset TP$ is a smooth distribution which is
complementary to the vertical subbundle, so we can consider $\om$ as
a generalized connection on $P$. Moreover $\om$ splits as
$\om=\om_{-1}+\om_{\goth b}$ according to the above decomposition and
similarly for the curvature.

A direct computation using property (2) of Cartan connections
shows that the curvature is always a
horizontal 2-form, i.e\. it vanishes if one of the vectors is
vertical. Thus it is fully described by the function
$\ka\in C^\infty(P, \goth
g^*_{-1}\otimes \goth g_{-1}^*\otimes \goth g)$,
$$
\ka(u)(X,Y)=K(\om^{-1}(X),\om^{-1}(Y))(u).
$$

If we evaluate the structure
equation on $\om^{-1}(X)$ and $\om^{-1}(Y)$ we obtain
$$\align
-[X,Y] + K(\om^{-1}(X), \om^{-1}(Y))
=\ &\om^{-1}(X)(\om(\om^{-1}(Y))) -\om^{-1}(Y)(\om(\om^{-1}(X)))
\\&-
\om([\om^{-1}(X), \om^{-1}(Y)])\\
=\ &- \om([\om^{-1}(X), \om^{-1}(Y)]).
\endalign$$
In particular for $X,Y\in \goth g_{-1}$, we see that
$\ka (u)(X,Y)=-\om(u)([\om^{-1}(X),\om^{-1}(Y)])$  so the $\goth b$-part of
$\ka$ is the obstruction against integrability of the horizontal
distribution defined by $\om$.

In particular, on the principal fiber bundle $G\to G/B$ over the
homogeneous space $G/B$ there is the (left) Maurer-Cartan form
$\om\in\Om^1(G,\goth g)$ which is a Cartan connection
and Maurer-Cartan structure equation
shows that the corresponding curvature is vanishing.

\subhead \nmb.{2.2}\endsubhead
Let $\la\: B\to GL(V_{\la})$ be a linear representation and let $P\to M$
be a principal fiber bundle as above. Then the
sections of the associated vector bundle $E_\la\to M$ correspond
bijectively to $B$-equivariant smooth functions in $C^\infty(P,
V_\la)$. We shall systematically use this identification without
further comments.
For a classical principal connection $\ga$ on $P$,
the covariant derivative on $E_\la$ along a vector field $X$ on $M$ can be
defined as the $B$-equivariant function $\nabla^\ga_X f$ with the
value at $u\in P$ given by the usual derivative in the direction
of the horizontal lift of $X$ to  $P$. For a general
connection on a bundle $P$ without any structure group, we can apply
the same definition to vector fields on $P$. The idea is to view the
Cartan connections as general connections,
but to exploit their special properties.

Given a Cartan connection $\om$ on $P$, there is the horizontal
projection $\chi_\om\:TP\to TP$ defined for each $\xi\in T_uP$
by $\chi_\om(\xi)= \xi - \ze_{\om_{\goth b}(\xi)}(u)$, where $\ze_Y$
means the fundamental vector field corresponding to $Y\in \goth b$. The
covariant exterior differential with respect to $\om$ on
vector-valued functions $s\in C^\infty(P,
V)$, evaluated on a vector field $\xi$ on $P$, is $d^{\om}
s(u)(\xi)=(\chi_\om\o \xi)\cdot s$. By definition, the value of the
covariant exterior differential $d^\om s(u)(\xi)$ with respect to
$\xi$ depends only on the horizontal projection $\chi_\om(\xi(u))$,
hence on $\om(\chi_\om(\xi(u)))=\om_{-1}(\xi(u))$. This leads to the
following

\subhead \nmb.{2.3}. Definition \endsubhead
The mapping $\nabla^{\om}\:C^\infty(P,V)\to C^\infty(P, \text{Hom}(\goth
g_{-1},V))$ defined by $
\nabla^\om s(u)(X) = \L_{\om^{-1}(X)} s(u) = (\om(u))^{-1}(X)\cdot s
$
is called the \idx{\it invariant differential} corresponding to the Cartan
connection $\om$.

We shall often use the brief notation $\nabla^\om_X s(u)$, $X\in
\goth g_{-1}$ for $(\nabla^\om s)(u)(X)$.

Note that the invariant differential of a $B$--equivariant function
is not $B$--equi\-va\-ri\-ant in general, so the invariant differential of
a section is not a section in general. Also, our brief notation
suggests that $\nabla^\om_X$ should behave like the usual covariant
derivative along a vector field, but the analogy fails in general
because of the nontrivial interaction between $\goth g_{-1}$ and
$\goth b$. There is however a possibility to form a section of a
bundle out of a given section and its invariant differential. This
point of view will be worked out in detail in section \nmb!{5}.

\proclaim{\nmb.{2.4}. Proposition (Bianchi identity)}
Let $\om\in\Om^1(P,\goth g)$ be a
Cartan connection.
Then the curvature $\ka$ satisfies
$$
\sum_{\text{cycl}}\bigl([\ka(X,Y),Z]-\ka(\ka_{-1}(X,Y),Z)-\nabla^\om_Z\ka(X,Y)
\bigr)=0
$$
for all $X$, $Y$, $Z\in\goth g_{-1}$, where $\sum_{\text{cycl}}$
denotes the sum over all cyclic permutations of the arguments.
\endproclaim
\demo{Proof} Let $X$, $Y$, $Z\in\goth g_{-1}$
and let us write $\tilde X$, $\tilde Y$, $\tilde Z$ for
the vector fields $\om^{-1}(X)$, $\om^{-1}(Y)$, $\om^{-1}(Z)$.
Now we evaluate the structure equation $d\om + \frac12[\om,\om]= K$
on the fields $[\tilde X,\tilde Y]$ and $\tilde Z$:
$$
-\L_{\tilde Z}\om([\tilde X,\tilde Y])-\om([[\tilde X,\tilde Y],\tilde Z])
- [\ka(X,Y),Z]
=- \ka(\ka_{-1}(X,Y),Z).
$$
Using the definition of the invariant differential, we obtain
$$
\om([[\tilde X,\tilde Y],\tilde Z]) =
-[\ka(X,Y),Z]+\ka(\ka_{-1}(X,Y),Z)+\nabla^\om_Z\ka(X,Y).
$$
Forming the cyclic sum, the left hand side vanishes by the Jacobi
identity for vector fields.\qed
\enddemo

\subhead \nmb.{2.5}. The iterated invariant differential \endsubhead
The invariant differential with respect to any Cartan connection $\om$
can be iterated, after $k$ applications on $s\in C^\infty(P,V)$
we get  $(\nabla)^ks = \nabla\dots\nabla s\in
C^\infty(P,\otimes ^k\goth g_{-1}^*\otimes V)$.
\proclaim{Lemma} For all $u\in P$ and $X$, $Y,\dots,Z\in\goth
g_{-1}$, $s\in C^\infty(P,V)$, we have
$$
(\nabla)^ks(u)(X,Y,\dots,Z) =
(\L_{\om^{-1}(Z)}\o\dots\o\L_{\om^{-1}(Y)}\o\L_{\om^{-1}(X)})s(u).
$$
In particular, we obtain
$(\nabla)^2s(u)(X,Y)-(\nabla)^2s(u)(Y,X)=\L_{\om^{-1}(\ka(X,Y))}s(u)$, the Ricci
identity.
\endproclaim
\demo{Proof} This is just the definition for $k=1$. So let us
assume that the statement holds for $k-1$. If we replace $s$ by its $(k-2)$-nd
invariant differential, we shall deal with the case $k=2$. By the
definition, $\nabla(\nabla s)(u)(X,Y)=\L_{\om^{-1}(Y)}(\nabla
s)(u)(X)=\L_{\om^{-1}(Y)}(\nabla s(\_)(X))(u)$ since the
invariant differential is linear in $X$. But the expression
in the last bracket is just $\L_{\om^{-1}(X)}s$.
Now,
$$
(\L_{\om^{-1}(Y)}\o\L_{\om^{-1}(X)}-
\L_{\om^{-1}(X)}\o\L_{\om^{-1}(Y)})s =
\L_{[\om^{-1}(Y),\om^{-1}(X)]}s =
\L_{\om^{-1}(\ka(X,Y))}s
$$
since $\goth g_{-1}$ is abelian.
\qed\enddemo

\subheading{\nmb.{2.6}} Let us compare our approach with the
classical covariant derivative with respect to a linear
connection $\ga\in\Om^1(P^1M,\goth g\goth l(m))$ on the linear
frame bundle $P^1M$. With the help of the soldering form $\th\in
\Om^1(P^1M,\Bbb R^m)$, we can build the form
$\om=\th\oplus\ga\in\Om^1(P^1M,\Bbb R^m\oplus \goth g\goth
l(m))$, a so called affine connection on $M$. It is simple to
check that $\om $ is a Cartan connection in the above sense and
that the horizontal lift of a vector $\xi\in T_xM$ which is
determined by $X\in \Bbb R^m$ and $u\in P^1M$ is exactly
$\om\i(X)\in T_u(P^1M)$. Thus the covariant differential
of a section $s$ of an associated bundle to $P^1M$ is given by
$\L_{\om\i(X)}\tilde s(u)$ where $\tilde s$ is the frame form of
$s$. Therefore the iterated differential $(\nabla^\om)^k$
coincides with the classical concept in this special case. The reason
why this case is much simpler than the general one is that $\Bbb R^m$
is an abelian ideal in $\Bbb R^m\oplus \goth g\goth l(m)$ and not
only a subalgebra.

\head \nmb0{3}. The second order structures
\endhead
\subhead \nmb.{3.1}\endsubhead
{}From now on we will assume that the group $G$ is connected and
semisimple, and that its Lie algebra is equipped with a grading
$\goth g=\goth g_{-1}\oplus \goth g_0\oplus \goth g_1$. Then the
following facts are well known, see \cite{Ochiai, 70}:
\roster
\item $\goth g_0$ is reductive with one--dimensional center
\item the map $\goth g_0\to\goth g\goth l(\goth g_{-1})$ induced by
       the adjoint representation is the inclusion of a subalgebra
\item the Killing form identifies $\goth g_1$ as a $\goth g_0$ module
       with the dual of $\goth g_{-1}$
\item the restrictions of the exponential map to $\goth g_1$ and
       $\goth g_{-1}$ are diffeomorphisms onto the corresponding
       closed subgroups of $G$.
\endroster

By $B$ we denote the closed (parabolic) subgroup of $G$ corresponding
to the Lie algebra $\goth b=\goth g_0\oplus \goth g_1$.
Then there is the normal subgroup $B_1$ in $B$ with Lie algebra
$\goth g_1$. From \therosteritem{4} above we see that $B_1$ is a
vector group. Finally, $B_0:=B/B_1$ is a reductive group with Lie
algebra $\goth g_0$, and the Lie group homomorphism induced by the
inclusion of $\goth g_0$ into $\goth b$ splits the projection, so $B$
is isomorphic to the semidirect product of $B_0$ and $B_1$.

\subheading{\nmb.{3.2}}
In this setting any Cartan connection
$\om\in \Om^1(P,\goth g_{-1}\oplus \goth g_0\oplus \goth g_1)$ on a
principal fiber bundle $P$ with structure group $B$ decomposes as
$\om = \om_{-1}\oplus\om_0\oplus \om_1$
and analogously does its curvature.

In order to involve certain covering phenomena, we shall slightly
extend  the classical definition of a
structure on a manifold $M$. For principal fiber bundles $P_1$,
$P_2$ over $M$ with structure groups $G_1$, $G_2$,
any morphism of principal fiber
bundles $P_1\to P_2$ over the identity on $M$, associated with
a covering of a subgroup of $G_2$ by $G_1$, will be
called a reduction of $P_2$ to the structure group $G_1$. For
example, the spin structures on Riemannian manifolds will be
incorporated into our general scheme in this way.

Now we show that in our setting, the canonical principal bundle
$G\to G/B=:M$ can be viewed as a reduction of the second order frame
bundle $P^2M\to M$ to the structure group $B$.

\proclaim{Lemma} Let $O\in M$ be the coset of $e\in G$,
$\ph\: \goth g_{-1}\to M$, $\ph(X)=\text{\rm exp}X\cdot O$. We define
$i\: G\to P^2M$, $i(g) = j^2_0(\ell_g\o \ph)$ and $i'\: B\to
G^2_m$, $i'(b) = j^2_0(\ph^{-1}\o \ell_b\o \ph)$. Then these two
mappings define a reduction (in the above sense) of $P^2M$.
\endproclaim
\demo{Proof} We have $i(g.b) = j^2_0(\ell_{g.b}\o \ph) =
j^2_0(\ell_g\o\ph).j^2_0(\ph^{-1}\o\ell_b\o\ph)$. Since the action of
$B$ on $M$ is induced by conjugation, the conditions
\nmb!{3.1}.\therosteritem{2} and \therosteritem{3}
imply that the homomorphism $i'$ induces an injection on the level
of Lie algebras, so it is indeed a covering of a subgroup of $G^2_m$.\qed
\enddemo
Before we give the general definition of a $B$--structure on a
manifold, we list some examples.

\subhead \nmb.{3.3}. Examples\endsubhead
The semisimple Lie algebras $\goth g$ which admit a grading of the form
$\goth g_{-1}\oplus \goth g_0\oplus \goth g_1$ can be completely
classified. In fact, the classification of these algebras in the
complex case is equivalent to the classification of Hermitian
symmetric spaces, see \cite{Baston, 91} for the relation. The full
classification can be found in \cite{Kobayashi, Nagano, 64, 65}
\sideremark{check the reference!}. Here we list some examples which
are of interest in geometry:

(1) Let $\goth g=\goth s\goth l(p+q, \Bbb R)$, the algebra of
matrices with trace zero, $\goth g_0=\goth s\goth l(p,\Bbb
R)\oplus \goth s\goth l(q,\Bbb R)\oplus \Bbb R$ and $\goth
g_{\pm1}= \Bbb R^{pq}$.  The grading is easily visible in a block form with
blocks of sizes $p,q$:
$$
\goth g_{-1} = \pmatrix 0&0\\*&0\endpmatrix,\newquad \goth
g_0=\pmatrix *&0\\0&*\endpmatrix,\newquad \goth g_1=\pmatrix
0&*\\0&0\endpmatrix.
$$
We obtain easily the formulae for the commutators. Let $X\in
\goth g_{-1}$, $Z\in \goth g_1$, $A=(A_1,A_2)\in\goth g_0$. Then
$$\alignat2
[~,~]&\:\goth g_0\x \goth g_{-1}\to \goth g_{-1},&\quad
[A,X]&=A_2.X-X.A_1\\
[~,~]&\:\goth g_1\x \goth g_{0}\to \goth g_{1},&\quad
[Z,A]&=Z.A_2-A_1.Z\\
[~,~]&\:\goth g_{-1}\x \goth g_{1}\to \goth g_{0},&\quad[X,Z]&=(-Z.X
,X.Z).
\endalignat$$

The corresponding homogeneous space is the real Grassmannian, the
corresponding structures are called \idx{\it almost Grassmannian}.
In the special case $p=1$, $q=m$, we obtain the classical
projective structures on $m$-dimensional manifolds.

(2) Let $\goth g=\goth s\goth o(m+1,n+1, \Bbb R)$, $\goth g_0 =
\goth c\goth o(m,n,\Bbb R)=\goth s\goth o(m,n,\Bbb R)\oplus\Bbb R$,
$\goth g_{-1}=\Bbb R^{m+n}$, $\goth
g_1=\Bbb R^{(m+n)*}$. For technical reasons we
choose the defining bilinear form
$\langle\ , \rangle$ on $\Bbb R^{m+n+2}$
given by $2x_0x_{m+n+1}+g(x_1,\dots,x_{m+n})$, where $g$ is the
standard pseudometric with signature $(m,n)$
given by the matrix $\Bbb J$. In block form with sizes $1,m+n,1$, we get
$$
\pmatrix 0&0&0\\p&0&0\\0&-p^T\Bbb J&0\endpmatrix\in \goth g_{-1}
,\newquad
\pmatrix -a&0&0\\0&A&0\\0&0&a\endpmatrix\in\goth g_0, \newquad
\pmatrix 0&q&0\\0&0&-\Bbb Jq^T\\0&0&0\endpmatrix\in\goth g_1,
$$
where $A\in \goth s\goth o(m,n,\Bbb R)$ and $a\Bbb I_{m+n}$ is in the center of
$\goth c\goth o(m,n,\Bbb R)$.

The commutators are
$$\alignat2
[~,~]&\:\goth g_{0}\x \goth g_{0}\to \goth
g_{0},&\quad[(A,a),(A',a')]&=(AA'-A'A,0)\\
[~,~]&\:\goth g_{0}\x \goth g_{-1}\to \goth
g_{-1},&\quad[(A,a),p]&=Ap+ap\\
[~,~]&\:\goth g_{1}\x \goth g_{0}\to \goth
g_{1},&\quad[q,(A,a)]&=qA+aq\\
[~,~]&\:\goth g_{-1}\x \goth g_{1}\to \goth
g_{0},&\quad[p,q]&=
(pq-\Bbb J(pq)^T\Bbb J,qp)
\endalignat$$
where $(A,a),\ (A',a')\in \goth s\goth o(m,n)\oplus\Bbb R=\goth g_0$,
$p\in\Bbb R^m= \goth g_{-1}$, $q\in\Bbb R^{m*}= \goth g_1$.

The homogeneous spaces are the
conformal pseudo--Riemannian spheres for metrics with signatures $(m,n)$.

(3) The symplectic algebra $\goth s\goth p(2n,
\Bbb R)$ admits the grading with $\goth g_{-1}=S^2\Bbb R^n$,
$\goth g_1=S^2\Bbb R^{n*}$, $\goth g_0=\goth g\goth l(n,\Bbb R)$.
We can express this grading in the block form:
$$
\pmatrix 0&X\\0&0\endpmatrix\in \goth g_{-1},\qquad
\pmatrix A&0\\0&-A^T\endpmatrix\in \goth g_0,\qquad
\pmatrix 0&0\\Z&0\endpmatrix\in \goth g_1
$$
The commutators are
$[X,Z]= X.Z\in \goth g\goth l(n,\Bbb R)$, $[A,X]=A.X+(A.X)^T\in
\goth g_{-1}$, $[A,Z]=-(Z.A+(Z.A)^T)\in \goth g_1$.
The corresponding homogenous spaces are the Lagrange Grassmann
manifolds and the corresponding structures are called almost
Lagrangian.

(4) If we use the symmetric form $\pmatrix 0 &\Bbb
I\\\Bbb I&0\endpmatrix$ instead of the antisymmetric one in the
previous example, then we obtain the grading $\goth s\goth
o(2n,\Bbb R)= \La^2\Bbb R^n\oplus \goth g\goth l(n,\Bbb R)\oplus
\La^2\Bbb R^{n*}$ with the commutators given by
$[X,Z]= X.Z\in \goth g\goth l(n,\Bbb R)$, $[A,X]=A.X-(A.X)^T\in
\goth g_{-1}$, $[A,Z]=-(Z.A-(Z.A)^T)\in \goth g_1$.
The corresponding homogeneous spaces are the isotropic Grassmann
manifolds. They can be identified with the spaces of pure spinors,
so the structures are called almost spinorial.

Some of the above examples coincide in small dimensions. Further
there are similar structures corresponding to the
exceptional Lie groups and we could also work in the
complex setting or choose different real forms.
For more information on these structures, see e.g\. \cite{Baston,
91}.

\subheading{\nmb.{3.4}. Definition}
Let $G$, $B$ be as in \nmb!{3.1} and let $M$ be a manifold of
dimension $m=\text{dim}(\goth g_{-1})$.  A {\it $B$-structure on $M$}
is a principal fiber bundle $P\to M$ with structure group $B$ which
is equipped with a differential form $\th =\th_{-1}\oplus
\th_0\in \Om^1(P,\goth g_{-1}\oplus\goth g_0)$ such that
\roster
\item $\th_{-1}(\xi)=0$ if and only if $\xi$ is a vertical
vector
\item $\th_0(\ze_{Y+Z}) = Y$ for all $Y\in \goth g_0$, $Z\in
\goth g_1$
\item $(r_b)^*\th = \operatorname{Ad}(b^{-1})\th$ for all $b\in B$,
where $\operatorname{Ad}$ means the action on the
vector space $\goth g_{-1}\oplus\goth g_0\simeq \goth g/\goth
g_1$ induced by the adjoint action.
\endroster
The form $\theta$ is called the {\it soldering form\/} or
{\it displacement form}. A homomorphism of $B$--structures is
just a homomorphism of principal bundles, which preserves the
soldering forms.

The {\it torsion\/} $T$ of the $B$--structure is defined by the structure
equation
$$
d\theta_{-1}=[\theta_{-1},\theta_0]+T
.$$

In the next part of this series, we shall apply the classical theory of
prolongations of $G$--structures to show that in all cases, except the
projective structures, each classical first order $B_0$--structure on $M$
gives rise to a distinguished $B$--structure on $M$ in the above sense. The
construction based on certain subtle normalizations extends essentially the
results on reductions of the second order frame bundles $P^2M$ due to
\cite{Ochiai, 70}, and it leads also to an explicit construction of the
canonical Cartan connection. An illustration of this procedure in the
special case of the conformal structures is presented in the last section of
this paper. The next lemma shows that in the Ochiai's approach we loose all
structures with non-zero torsion $T$. However, the torsion is quite often
the only obstruction against the local flatness, see \cite{Baston, 91} or
\cite{\v Cap, Slov\'ak, 95}. For example, Ochiai deals in fact only with
spaces locally isomorphic to the
homogeneous spaces for all higher dimensional Grassmannian structures.

\proclaim{\nmb.{3.5}. Lemma}
Let $M$ be an $m$--dimensional manifold and let $P\to M$ be a
reduction of the second order frame bundle $P^2M\to M$ to the group
$B$ over the homomorphism $i'$, as in \nmb!{3.2}. Then there is a
canonical soldering form $\theta$ on $P$ such that $(P,\theta)$ is
a $B$--structure on $M$, and this $B$--structure has zero torsion.
\endproclaim
\demo{Proof}
The second frame bundle $P^2M$, is equipped with a canonical
soldering form a form
$\th^{(2)}\in\Om^1(P^2M,\Bbb R^m\oplus\goth g^1_m)$ defined as
follows. Each element
$u\in P^2_xM$, $u=j^2_0\ph$, determines a linear isomorphism
$\tilde u\:\Bbb R^m\oplus \goth g^1_m\to T_{\pi^2_1(u)}P^1M$ (in fact
$T_0(P^1\ph)\:T_{(0,e)}(\Bbb R^m\x G^1_m)\to TP^1M$). Now if
$X\in T_uP^2M$
then $\th^{(2)}(X)=\tilde u\i(T\pi^2_1(X))$, i.e\.
$\th^{(2)}(X)=j^1_0(P^1\ph\i\o\pi^2_1\o c)$ if $X=j^1_0c$. This
canonical form
decomposes as $\th^{(2)}=\th_{-1}\oplus\th_0$ where $\th_{-1}$ is the
pullback of the soldering form $\th$ on $P^1M$,
$\th_{-1}=(\pi^2_1)^*\th$, while $\th_0$ is $\goth g^1_m$-valued.

It is well known, see e.g\. in \cite{Kobayashi, 72} that this is in fact a
soldering form with zero torsion, and that for any
reduction as assumed the pullback of $\theta$ is again a
soldering form, which clearly has trivial torsion, too.
\qed\enddemo

\subhead \nmb.{3.6} \endsubhead
A $B$--structure $(P,\theta)$ is related to a rich underlying
structure. First, we can form the bundle $P_0:=P/B_1\to M$, which is
clearly a principal bundle with group $B_0$, and $P\to P/B_1$ is a
principal $B_1$--bundle. Now consider the component
$\theta_{-1}$ of the soldering form. By property \therosteritem{3}
in the definition of the soldering form it is $B_1$--invariant and
clearly it is horizontal as a form on $P\to P_0$, so it passes down
to a well defined form in $\Om^1(P_0,\goth g_{-1})$, which we again
denote by $\theta_{-1}$. One easily verifies that this form is
$B_0$--equivariant and its kernel on each tangent space is precisely
the vertical tangent space of $P_0\to M$. Then for each $u\in P_0$,
$\theta_{-1}$ induces a linear isomorphism
$T_uP_0/V_uP_0\simeq \goth g_{-1}$ and composing the inverse of this
map with the tangent map of the projection $p:P_0\to M$, we associate
to each $u\in P_0$ a linear isomorphism
$\goth g_{-1}\simeq T_{p(u)}M$, thus obtaining a reduction
$P_0\to P^1M$ of the frame bundle of $M$ to the group $B_0$, where
$B_0$ is mapped to $\text{GL}(m,\Bbb R)\simeq \text{GL}(\goth g_{-1})$
via the adjoint action.

In particular this shows that one can view the tangent bundle $TM$ of
$M$ as the associated bundle $P_0\x_{\text{\rm Ad}}\goth g_{-1}$. Since
$P_0=P/B_1$, we can as well identify $TM$ with $P\x_{(\text{\rm Ad},
\text{\rm id})}\goth g_{-1}$. Here (Ad,id)  means the adjoint action
of $B_0$ and the trivial action of $B_1$.

\proclaim{Lemma} Let $(P,\theta)$ be a $B$--structure on $M$,
$P_0:=P/B_1$. Then there exists a global smooth $B_0$--equivariant
section $P_0\to P$, and if $\si$ is any such section we have:
\roster
\item $\ga := \si^*\theta_0\in \Om^1(P_0,\goth g_0)$ is a
       principal connection on $P_0$.
\item $\om:=\si^*\theta$ is a Cartan connection on $P_0$ with
       $\goth g_{-1}$--component equal to the form $\theta_{-1}$ from above.
\item The invariant differential $\nabla^{\om}\: C^\infty(P_0,V)\to
       C^\infty(P_0, \goth g_{-1}^*\otimes V)$ coincides with the usual
       covariant (exterior) differential
       $d^{\ga}\:\Om^0(P_0,V)\to \Om^1(P_0,V)$ viewed as
       $d^\ga\:C^\infty(P_0,V)^{B_0}\to C^\infty(P_0, \goth g_{-1}^*\otimes
       V)^{B_0}$.
\item The components of the curvature $K=K_{-1}\oplus K_0$ of $\om$
       are just the torsion and the curvature of the principal
       connection $\ga$.
\endroster
The space of all equivariant sections $\si$ as above is an affine
space modeled on the space $\Om^1(M)$ of one--forms on $M$.
\endproclaim
\demo{Proof}
Starting from a principal bundle atlas for $P\to M$, we see that we
can find a covering $\{U_\al\}$ of $M$ such that the bundle $p:P\to
P_0$ is trivial over any of the sets $\pi^{-1}(U_\al)\subset P_0$,
where $\pi:P_0\to M$ is the projection. Since $B$ is the semidirect
product of $B_0$ and $B_1$ we can choose a local $B_0$--equivariant
section $s_\al$ of $P\to P_0$ over each of these subsets.

Next, there is a smooth mapping $\ch:P\x _{P_0}P\to \goth g_1$
determined by the equation $v=u\cdot\text{exp}(\ch(u,v))$. If
$U_\al\cap U_\be\neq \emptyset$ then we have a well defined smooth map
$\ch_{\al\be}:\pi^{-1}(U_\al\cap U_\be)\to \goth g_1$ given by
$\ch_{\al\be}(u):=\ch(s_\al(u),s_\be(u))$. Since the sections are
$B_0$--equivariant one easily verifies that $\ch_{\al\be}(u\cdot
b)=\Ad(b^{-1})\cdot \ch_{\al\be}(u)$ for all $b\in B_0$. Let
$\{f_\al\}$ be a partition of unity subordinate to the covering
$\{U_\al\}$ of $M$. For $u\in P_0$ define $s(u)\in P$ as follows:
Choose an $\al$ with $\pi(u)\in U_\al$ and put
$$
s(u):=s_\al(u)\cdot\text{exp}(\tsum_{\be}f_\be(\pi(u))\ch_{\al\be}(u)).
$$
Clearly this expression makes sense, although the $\ch_{\al\be}$ are
only locally defined. Since $B_1$ is abelian, it is easily seen that
$\ch_{\al\ga}(u)=\ch_{\al\be}(u)+\ch_{\be\ga}(u)$, whenever all terms
are defined. Now if $\ga$ is another index such that
$\pi(u)\in U_\ga$, we get:
$$\multline
s_\ga(u)\cdot\text{exp}(\tsum_{\be}f_\be(\pi(u))\ch_{\ga\be}(u))=\\
=s_\al(u)\cdot\text{exp}(\ch_{\al\ga}(u))\cdot
\text{exp}(\tsum_{\be}f_\be(\pi(u))(\ch_{\ga\al}(u)+\ch_{\al\be}(u)))=\\
=s_\al(u)\cdot\text{exp}(\ch_{\al\ga}(u)+\ch_{\ga\al}(u)+
\tsum_{\be}f_\be(\pi(u))\ch_{\al\be}(u)),
\endmultline$$
so $s(u)$ is independent of the choice of $\al$, and thus $s:P_0\to P$
is a well defined smooth global section. Moreover, for $b\in B_0$
$$
s(u\cdot b)=s_{\al}(u)\cdot b\cdot\text{exp}(\Ad(b)\cdot
\tsum_{\be}f_\be(\pi(u))\ch_{\al\be}(u))=s(u)\cdot b,
$$
so $s$ is $B_0$--equivariant, too.

Now if $s$ and $\si$ are two global equivariant sections, then
$u\mapsto \ch(s(u),\si(u))$ is the frame form of a smooth one--form on
$M$. On the other hand if $\ph:P_0\to \goth g_1$ is the frame form of
a one--form, then $u\mapsto s(u)\cdot \operatorname{exp}(\ph(u))$
is again a smooth equivariant section.

\therosteritem{1} and \therosteritem{2} are easily verified directly,
and \therosteritem{3} was shown in \nmb!{2.6}. \therosteritem{4}
follows immediately from \therosteritem{3} and \therosteritem{2}
since $\om=\theta_{-1}\oplus \gamma$ and the torsion and curvature
of $\ga$ are by definition just $d^\ga\theta_{-1}$ and $d^\ga\ga$,
respectively.
\qed\enddemo

\subhead \nmb.{3.7}. Induced Cartan connections \endsubhead
We have seen in \nmb!{3.6} that the soldering form on $P\to M$
leads via $B_0$--equivariant sections $\si:P_0\to P$ to a
distinguished class of principal connections on the bundle
$P_0\to M$, which can be canonically extended to Cartan connections
on the latter bundle. Next we show that to any principal connection
$\ga$ from this distinguished class, i.e\. to each equivariant
section $\si$ as above, we can construct an induced Cartan connection
$\tilde \ga$ on $P$, which is $\si$--related to the Cartan connection
$\theta_{-1}\oplus \ga$.

\proclaim{Lemma}
For each $B_0$-equivariant section $\si\:P_0\to P$, there is
a uniquely defined Cartan connection
$\om=\theta_{-1}\oplus\theta_0\oplus\om_1$ satisfying
$\om_1|(T\si(TP_0))=0$.
\endproclaim
\demo{Proof}
Using condition \therosteritem{4} of \nmb!{3.1} we see that
the section $\si$ induces an isomorphism of $P$
with $P_0\x \goth g_1$ defined by $u\mapsto (p(u), \ta(u))$ where
$p:P\to P_0$ is the projection and the mapping
$\ta\: P\to \goth g_1$ is defined by the equality
$u = \si(p(u)).\text{exp}(\ta(u))$.
Now we define $\om_1$ on $\si (P_0)$ by $\om_1|_{\si(P_0)}:=d\tau$.
Since $\tau\o \si=0$ we clearly have $\om_1|(T\si(TP_0))=0$,
and obviously for any $u\in \si(P_0)$, $\om$ induces a
bijection $T_uP\to \goth g$.

Next, since $\tau$ is identically zero on $\si(P_0)$ and
$\tau(u\cdot \text{\rm exp}tX)=\tau(u)+tX$ for
$X\in \goth g_1$ it follows from \nmb!{3.4}.\therosteritem{2} that
on $\si(P_0)$ the form $\om=\theta\oplus \om_1$ reproduces the
generators of fundamental fields.

Now one easily checks that this form can be uniquely extended using
the equivariancy properties which are required for a Cartan
connection.
\qed\enddemo

\proclaim{\nmb.{3.8}. Lemma} In the situation of \nmb!{3.7}, denote
by $p:P\to P_0$ the projection and let $V$ be any representation of
$B_0$. Let $u\in P$ and $X,\,Y\in\goth g_{-1}$,
$b=\text{\rm exp}(\ta(u))$, where $\tau$ is the mapping from the
proof of \nmb!{3.7}, and $s\in C^\infty(P_0,V)$. Then we have:
\roster\item $(\nabla^{\tilde \ga}\o p^* - p^*\o\nabla^{\ga})s
(u)(X) = \ze_{[\ta(u),X]}(\si(p(u))).(s\o p)$.
\item Let $h\in B$ be arbitrary.
The curvature $\ka\in C^\infty(P,\goth g_{-1}^*\otimes \goth
g_{-1}^*\otimes \goth g)$ of any Cartan connection satisfies
$\ka(X,Y)(u.h)=\Ad(h^{-1}).\ka(\Ad (h).X, \Ad
(h).Y)(u)$.
\item The curvature $K$ of $\tilde\ga$ and the curvature
$R$ of $\ga$ satisfy $d^\ga\th_{-1}\oplus R=\si^*K$. In particular
$R=\si^* K_0$.
\item The curvature components of an arbitrary Cartan connection satisfy
$$\align
\ka_{-1}(u)(X,Y)&=\ka_{-1}(\si(p(u)))(X,Y)\\
\ka_0(u)(X,Y) &=
\ka_0(\si(p(u)))(X,Y)-[\ta(u),\ka_{-1}(\si(p(u)))(X,Y)]\\
\ka_1(u)(X,Y) &= \ka_{1}(\si(p(u)))(X,Y)-[\ta(u), \ka_0(\si(p(u)))(X,Y)]
+\\&\newquad\frac12[\ta(u),[\ta(u),\ka_{-1}(\si(p(u)))(X,Y)]].
\endalign$$
\item If the $B$--structure has zero torsion, then the curvature of
the induced Cartan connection $\tilde \ga$ satisfies
$$\align
\ka_{-1}(u)(X,Y)&=0\\
\ka_0(u)(X,Y) &= \ka_0(\si(p(u)))(X,Y)\\
\ka_1(u)(X,Y) &= [\ka_0(\si(p(u)))(X,Y),\ta(u)].
\endalign$$
In particular, the component $\ka_1$ vanishes on $\si(P_{0})$.
\endroster
\endproclaim

\demo{Proof} By the definition of the Cartan connections, the
formula
$$
Tr^b(\tilde \ga^{-1}(X)(u))=
\tilde \ga^{-1}(\Ad(b^{-1}).X)(u.b)
\tag6$$
holds for all $u\in P$, $b\in B$.
Since $\ga = \si^*(\tilde\ga)_0$, the horizontal lift of the
vector $\xi\in T_xM$ corresponding to $X\in \goth g_{-1}$ and
$p(u)\in P_0$ with respect to $\ga$ is just $Tp(\tilde
\ga\i(X)(\si(p(u))))$, see \nmb!{2.6}.
The definition of the
tangent mapping then yields
$$\align
(\nabla^{\tilde \ga}\o p^* &- p^*\o\nabla^{\ga})s
(u)(X) = \tilde\ga^{-1}(X)(u).(s\o p) - Tp\bigl(\tilde
\ga^{-1}(X)(\si\o p(u))\bigr). s\\
&= Tr^b(\tilde\ga^{-1}(\Ad b.X)(\si\o p(u))).(s\o p) - Tp\bigl(
\tilde\ga^{-1}(X)(\si\o p(u))\bigr).s\\
&= \tilde \ga^{-1}(X +[\ta(u),X])(\si\o p(u)).(s\o p) - Tp\bigl(
\tilde\ga^{-1}(X)(\si\o p(u))\bigr).s\\
&= \ze_{[\ta(u),X]}(\si\o p(u)).(s\o p)
\endalign$$
where the last but one equality is obtained using the fact that for
$Z\in\goth g_1$, $X\in\goth g_{-1}$ we have:
$$\aligned
(\text{Ad}(\text{exp}Z)).X &= X + [Z,X] + \frac12[Z,[Z,X]] +
\frac16[Z,[Z,[Z,X]]]+\dots\\
&= X + [Z,X] + \frac12[Z,[Z,X]].
\endaligned\tag 7$$

The next claim also follows from the formula (6) and from the fact
that the Lie bracket of $f$-related vector fields is $f$-related:
$$\align
\ka(X,Y)(u.b) &= K(\tilde\ga^{-1}(X),\tilde \ga^{-1}(Y))(u.b)\\
&=-\tilde\ga([\tilde\ga^{-1}(X),\tilde\ga^{-1}(Y)](u.b))\\
&= -\Ad(b^{-1})\o \tilde\ga(Tr^{b^{-1}}.
([\tilde\ga^{-1}(X),\tilde\ga^{-1}(Y)](u.b)))\\
&= -\Ad(b^{-1})\o\tilde\ga([\tilde\ga^{-1}(\Ad
b.X),\tilde\ga^{-1}(\Ad b.Y)](u))\\
&=\Ad (b_{-1}).\ka(\Ad b.X,\Ad b.Y)(u).
\endalign$$

(3) follows immediately from the fact that $\tilde \ga$ and
$\theta_{-1}\oplus \ga$ are $\si$--related.

If $b\in \goth g_1$, the horizontal part of $\tilde \ga^{-1}(\Ad
b.X)$ is just $\tilde \ga^{-1}(X)$. The curvature of a Cartan
connection is a horizontal form and so (2) implies that
$\ka(X,Y)(u)=\Ad b^{-1}\ka(X,Y)(\si(p(u)))$.
Now \nmb!{3.8}.(7) implies directly the relations (4).

Once we prove that $\ka_1|\si(P_0)=0$, (5) will follow
directly form (4) since in this case $\ka_{-1}$ is just the torsion.
But according to the definition of $\tilde \ga$, the vector fields
$\tilde\ga^{-1}(X)$ are tangent to $\si(P_0)$ for all $X\in \goth g_{-1}$.
Consequently also the Lie brackets of such fields are tangent to $\si(P_0)$
and thus $\tilde \ga_1([\tilde\ga^{-1}(X),\tilde\ga^{-1}(Y)])=0$.
\qed\enddemo

\subhead \nmb.{3.9}. Admissible Cartan connections \endsubhead
Let $(P,\theta)$ be a $B$--structure on $M$. A Cartan connection
$\om$ on $P$ is called {\it admissible\/} if and only if it is of the
form $\om =\th_{-1}\oplus\th_0\oplus\om_1$. Thus in particular the
induced connections from \nmb!{3.7} are admissible. Moreover, by
definition the $\goth g_{-1}$ component of the curvature of any
admissible Cartan connection is given by the torsion of the
$B$--structure.

Let us now consider two admissible Cartan connections $\om$, $\bar
\om$, so that they differ only in the $\goth g_1$-component. Then
there is a function
$\Ga\in C^\infty(P,\goth g_{-1}^*\otimes \goth g_{1})$
such that $\bar \om=\om - \Ga\o \th_{-1}$. Indeed, $\om - \bar\om$ has
values in $\goth g_1$ and vanishes on vertical vectors.

The function $\Ga$ can be viewed as an expression for the
``deformation'' of $\om$ into $\bar\om$ and in view of its
properties proved below, we call it the {\it deformation tensor}.

\proclaim{\nmb.{3.10}. Lemma} \roster\runinitem $\Ga(u.b) = \text{\rm
Ad}(b^{-1})\o \Ga(u)\o \text{\rm Ad}(b)$ for all $b\in B_0$
\item $\Ga(u.b)=\Ga(u)$ for all $b\in B_1$
\item $\bar \om^{-1}(X)(u) = \om^{-1}(X)(u) + \ze_{\Ga(u).X}(u)$ for
all $X\in \goth g_{-1}$
\item $(\ka_0-\bar \ka_0)(u)(X,Y)= [X,\Ga(u).Y]+[\Ga(u).X,Y]$
\item $(\ka_1-\bar\ka_1)(u)(X,Y)=
\nabla^{\om}_X\Ga(u).Y-\nabla^{\om}_Y\Ga(u).X +
\Ga(u)(\ka_{-1}(X,Y))$
\item $(\ka_{-1}-\bar\ka_{-1})(u)(X,Y)=0$.
\endroster
\endproclaim
\demo{Proof}
By definition,
$(r^b)^*(\Ga\o\th_{-1})=\text{Ad}(b^{-1})\o(\Ga\o\th_{-1})$
and the adjoint action is trivial if $b\in B_1$. Since
$(r^b)^*\th_{-1} = \th_{-1}$ for $b\in B_1$ too,
the second claim  has been proved.
If $b\in B_0$ then $\Ga\o\th_{-1}(Tr^b.\xi)(u.b)
=\text{Ad}(b^{-1})(\Ga\o\th_{-1})(\xi)(u)$ and the left hand side is
$\Ga(u.b)\o\text{Ad}(b^{-1})\o\th_{-1}(u)(\xi)$ by the equivariancy of
$\th_{-1}$. Comparing the results, we obtain just the required
formula (1).

In order to obtain (3), we compute $\bar \om(\om^{-1}(X)) = X
- \Ga\o\th_{-1}(\om^{-1}(X))$ and so $\bar \om^{-1}(X)=
\om^{-1}(X)+\bar
\om^{-1}(\Ga\o\th(\om^{-1}(X)))=\om^{-1}(X)+\ze_{\Ga.X}$.

In order to verify (4) and (5), let us compute (we use just the
definition of the frame form of the curvature)
$$
\align
(\ka&-\bar\ka)(X,Y) =
\bar\om([\bar\om^{-1}(X),\bar\om^{-1}(Y)])-
\om([\om^{-1}(X),\om^{-1}(Y)])\\
&= (\om -
\Ga\o\th)([\om^{-1}(X)+\ze_{\Ga.X},\om^{-1}(Y)+\ze_{\Ga.Y}]) -
\om([\om^{-1}(X),\om^{-1}(Y)])\\
&= \om([\ze_{\Ga.X},\om^{-1}(Y)]) + \om([\om^{-1}(X),\ze_{\Ga.Y}]) +
\om([\ze_{\Ga.X},\ze_{\Ga.Y}])+\Ga(\ka_{-1}(X,Y))-\\
&\newquad \Ga\o\om_{-1}([\ze_{\Ga.X},\om^{-1}(Y)]) -
\Ga\o\om_{-1}([\om^{-1}(X),\ze_{\Ga.Y}]) -
\Ga\o\om_{-1}([\ze_{\Ga.X},\ze_{\Ga.Y}]).
\endalign$$
We have to notice that the fields
$\ze_{\Ga.X}(u)=\om^{-1}(\Ga(u).(X))(u)$ are defined by
means of the fundamental field mapping, but with arguments
varying from point to point in $P$.
To resolve the individual brackets, we shall evaluate the
curvature $K$ of $\om$ on the corresponding fields:
$$
\align
d\om&(\om^{-1}(\Ga.X),\om^{-1}(Y)) =\\
&=\L_{\om^{-1}(\Ga.X)}\om(\om^{-1}(Y)) -
\L_{\om^{-1}(Y)}\om(\om^{-1}(\Ga(X)))-
\om([\om^{-1}(\Ga(X)),\om^{-1}(Y)])\\
&=-[\om(\om^{-1}(\Ga(X))),\om(\om^{-1}(Y))] +
K(\om^{-1}(\Ga(X)),\om^{-1}(Y)).
\endalign$$
Since $K$ is a horizontal 2-form it evaluates to zero and
$\om(\om^{-1}(Y))=Y$ is constant. Thus we obtain
$$
\om([\om^{-1}(\Ga(X)),\om^{-1}(Y)]) = [\Ga.X,Y]-
\L_{\om^{-1}(Y)}\Ga.X.$$
Now we can decompose this equality into the individual
components.
$$\align
\om_{-1}([\om^{-1}(\Ga(X)),\om^{-1}(Y)] &= 0\\
\om_{0}([\om^{-1}(\Ga(X)),\om^{-1}(Y)]  &= [\Ga.X,Y]\\
\om_{1}([\om^{-1}(\Ga(X)),\om^{-1}(Y)]  &= -\nabla^{\om}_Y\Ga.X.
\endalign$$
It remains to evaluate the structure equation on the fields
$\ze_{\Ga.X}$, $\ze_{\Ga.Y}$. Since $\goth g_1$ is abelian and
$K(\ze_{\Ga.X}, \ze_{\Ga.Y})=0$, we obtain
$$
\om([\om^{-1}(\Ga(X)),\om^{-1}(\Ga(Y))])=
\L_{\om^{-1}(\Ga(X))}\Ga.Y -
\L_{\om^{-1}(\Ga(Y))}\Ga.X.
$$
But the Lie derivatives of $\Ga$  depend only on the value of the
vector field in the point in question. However, we have already
proved that $\Ga$ is $B_1$-invariant and consequently the
derivative is zero. Now we can insert the expressions for the
brackets into the above expression for the difference $\ka-\bar\ka$
and we get exactly the required formulae.
\qed
\enddemo

In particular (1) and (2) show that
$\Ga$ is always a pullback of a tensor on $M$. This fact is of
basic importance for our approach.

\head \nmb0{4}. Formulae for the iterated invariant differential
\endhead

As before, we shall consider sections $s\in C^\infty(P_0, V_\la)^{B_0}$ of
associated bundles induced by representations of $B_0$ and we shall view
them as equivariant mappings $p^*s\in C^\infty (P,V_\la)^B$.  We shall
develop a recurrence procedure which expands the iterated differentials of
such sections with respect to any admissible connection in terms of the
underlying linear connections. This expression splits the invariant
derivatives into the equivariant part (thus a section) and the obstruction
parts (which concentrates the failure to the $B_1$-invariance).
Thus, after having the canonical Cartan connections,
this will provide us with a direct method of constructing the invariant
operators.

\subheading{\nmb.{4.1}}
In view of the results of the preceding section, the comparison of
the iterated covariant
differential with respect to the principal connection
$\ga=\si^*\om_0$ on $P_0$, with the invariant differential
$\nabla^\om$ with respect to the admissible
Cartan connection $\om$ becomes quite algorithmic. Indeed
we can write
$$\align
(\nabla^{\om})^k\o p^* &-p^* \o (\nabla^{\ga})^k =
\nabla^{\om}\o((\nabla^{\om})^{k-1}\o p^* -
p^* \o(\nabla^\ga)^{k-1})+\\&\newquad
(\nabla^{\om}\o p^* - p^*\o \nabla^\ga)\o(\nabla^\ga)^{k-1}\\
&=\nabla^{\om}\o((\nabla^{\om})^{k-1}\o p^* -
p^* \o(\nabla^\ga)^{k-1})+
\\&\newquad
(\nabla^{\om}\o p^* - \nabla^{\tilde\ga}\o
p^*)\o(\nabla^\ga)^{k-1}+(\nabla^{\tilde\ga}\o p^* -
p^*\o\nabla^{\ga})(\nabla^\ga)^{k-1}.
\endalign$$

Thus, we have to start an induction procedure.
Let us remind, that the deformation tensor $\Ga\in \C^\infty(P,\goth
g_{-1}^*\otimes \goth g_1)$ transforming $\tilde \ga$
into $\om$ is a pullback of a tensor on $M$, see Lemma \nmb!{3.10}.
We shall work in the setting of \nmb!{3.6}-\nmb!{3.10} with $s\in
C^\infty(P_0, V_{\la})^{B_0}$, where $V_{\la}$ is the representation space
for $\la\:B_0\to GL(V_{\la})$.
In particular, we know from Lemmas \nmb!{3.8}, \nmb!{3.10} that for
$u\in P$, $X,\,Y\in\goth g_{-1}$
$$\gather
\nabla^\om_X(p^*s)(u)-\nabla^{\tilde\ga}_X(p^*s)(u) =
\ze_{\Ga(u).X}(u).(p^*s)=0\\
(\nabla^{\tilde \ga}\o p^* - p^*\o\nabla^{\ga})s
(u)(X) = \ze_{[\ta(u),X]}(\si(p(u))).(p^*s)=\la([X,\ta(u)])(s(p(u))).
\endgather$$
Consequently, the middle term in the above inductive formula vanishes and
the last one yields always the induced action of the bracket on
the target space of the iterated covariant differential
$(\nabla^\ga)^{k-1}$. In particular,
we have already deduced the general formula for the first order
operators:

\proclaim{\nmb.{4.2}. Proposition} Let $\om$ be an admissible Cartan
connection, $\ga$ be the linear connection corresponding to
an equivariant section $\si\:P_0\to P$. For all $X\in \goth
g_{-1}$, $s\in C^{\infty}(P_0, V_\la)^{B_0}$, $u\in P$ we have
$$\
(\nabla^\om\o p^* - p^*\o \nabla^\ga)s(u)(X)=
\la([X,\ta(u)])(s(p(u))).
$$
In particular, the difference is zero if evaluated at points with
$\ta(u)=0$.
\endproclaim

In order to continue to higher orders, we need to know
how to differentiate the expressions which will appear. Thus let
us continue with two technical lemmas.

\proclaim{\nmb.{4.3}. Lemma} Let $X\in \goth g_{-1}$,
$Z\in\goth g_1$, $u\in P$. Then
$$\align
\L_{\tilde\ga^{-1}(X)}\ta(u)&= \tfrac12[\ta(u),[\ta(u),X]]\tag1\\
\L_{\om^{-1}(Z)}\ta(u) &= Z\tag2
\endalign
$$
\endproclaim
\demo{Proof} The definition of $\ta$ can be written as
$\ta(\si(p(u))=0$, $\ta\o r^{\text{exp}Z}(u) =\ta(u)+Z$, $Z\in
\goth g_1$. Thus in order to get (1), we can compute for $u=\si(p(u)).b$
$$\align
\nabla^{\tilde\ga}_X&\ta(u)=T\ta.(\tilde\ga^{-1}(X)(u))=
T\ta\o Tr^b.\tilde\ga^{-1}(\text{Ad}b.X)(\si(p(u)))\\&=
T(\ta\o r^b)(\tilde\ga^{-1}(X + [\ta(u), X] + \tfrac12[\ta(u),
[\ta(u), X]])(\si(p(u))))\\&=
T\ta(\tilde\ga^{-1}(X + [\ta(u), X])(\si(p(u)))) + \tfrac12T\ta(\tilde\ga^{-1}([\ta(u),
[\ta(u), X]])(\si(p(u))))\\&=
\tfrac12[\ta(u),[\ta(u), X]].
\endalign$$
Next let us compute $\ze_{Z}.\ta(u)$  for $Z\in \goth g_1$.
$$
T\ta(\ze_{Z})(u)= \left.\tfrac\partial{\partial t}\right|_0\bigl(
\ta(u.\text{exp}tZ\bigr)
= \left.\tfrac\partial{\partial t}\right|_0\bigl(
\ta(u)+tZ\bigr) = Z
$$
i.e\. (2) holds.
\qed\enddemo

\proclaim{\nmb.{4.4}. Lemma} Let $f\:P\to V_{\la}$ be a
mapping defined by
$$
f(u) = \tilde f(p(u))(\ta(u),\dots, \ta(u)),
$$
where $\tilde f\:P_0\to \otimes^k\goth g_1^*\otimes V_\la$ is
$\goth g_0$-equivariant with respect to the canonical action
$\tilde \la$ on the tensor product. Then
$$\align
\nabla^{\om}_Y f(u) &= \la([Y,\ta(u)])(f(u)) - \\
&\newquad\tfrac12\tsum_{i=1}^k(p^*\tilde
f)(u)(\ta(u),\dots,[\ta(u),[\ta(u),Y]],\dots,\ta(u)) +\\
&\newquad (p^*(\nabla^\ga_Y \tilde f))(u)(\ta(u),\dots,\ta(u))
+\\&\newquad
\tsum_{i=1}^k (p^*\tilde f)(u)(\ta(u),\dots,\Ga(u).Y,\dots,\ta(u)).
\endalign$$
Moreover, all the terms in the above expression for $\nabla^\om
f\:P\to \goth g_{-1}^*\otimes V_\la$ satisfy the assumptions of this lemma with the corresponding
canonical representation on $\otimes^t\goth g_1^*\otimes\goth
g_{-1}^*\otimes V_\la$, where $t$ is the number of $\ta$'s entering
the term in question. \endproclaim

\demo{Proof} Let us compute using the chain rule, Proposition
\nmb!{4.2}, Lemma \nmb!{4.3} and \nmb!{4.1}
$$\align
(\nabla^\om_Y f)(u) &= (\nabla^\om_Y(p^*\tilde
f)(u))(\ta(u),\dots,\ta(u)) +\\
&\newquad \tsum_{i=1}^k(p^*\tilde
f)(u)(\ta(u),\dots,\nabla^\om_Y\ta(u),\dots,\ta(u))\\
&=(p^*(\nabla^\ga_Y\tilde f)(u))(\ta(u),\dots,\ta(u)) +\\
&\newquad (\tilde
\la([Y,\ta(u)])(p^*\tilde f)(u))(\ta(u),\dots,\ta(u)) +\\
&\newquad \tfrac12\tsum_{i=1}^k (p^*\tilde
f)(u)(\ta(u),\dots,[\ta(u),[\ta(u),Y]],\dots, \ta(u)) + \\
&\newquad \tsum_{i=1}^k (p^*\tilde
f)(u)(\ta(u),\dots,\Ga(u).Y,\dots, \ta(u))\\
&=p^*(\nabla^\ga_Y\tilde f) (u)(\ta(u),\dots,\ta(u)) + \\
&\newquad \la([Y,\ta(u)])(f(u)) -\\
&\newquad \tsum_{i=1}^k(p^*\tilde
f)(u)(\ta(u),\dots,[[Y,\ta(u)],\ta(u)],\dots,\ta(u)) +\\
&\newquad \tfrac12\tsum_{i=1}^k (p^*\tilde
f)(u)(\ta(u),\dots,[[Y,\ta(u)],\ta(u)],\dots, \ta(u)) +\\
&\newquad\tsum_{i=1}^k (p^*\tilde
f)(u)(\ta(u),\dots,\Ga(u).Y,\dots, \ta(u)).
\endalign$$
It remains to prove that the resulting expressions satisfy once
more the assumptions of the lemma. Let us show the argument on
the first term $f_1(u)(Y) := \la([Y,\ta(u)])f(u)$. We have
$f_1(u)(Y) = \tilde f_1(p(u))(\ta(u),\dots,\ta(u))(Y)$ with $\tilde
f_1\:P_0\to \otimes^{k+1}\goth g_1^*\otimes \goth g_{-1}^*\otimes
V_\la$, $\tilde f_1(p(u))(Z_0,\dots,Z_k,Y) = \la([Y,Z_0])(\tilde
f(p(u))(Z_1,\dots,Z_k))$. The evaluation of $\tilde f_1$ on $Z_0$
and $Y$ can be written as the composition
$$
(\text{id}\otimes\la)\o(\tilde f\otimes \text{ad})\:P_0\x \goth
g_1\otimes\goth g_{-1}\to \otimes^{k}\goth g_1^*\otimes V_\la
$$
of equivariant mappings, so $\tilde f_1$ is equivariant as
well. Similarly one can write down explicitly the terms in the
second and the third part of the expression. The equivariancy of
the terms with $\Ga$ follows from \nmb!{3.10}.
\qed\enddemo

\subheading{\nmb.{4.5}. The second order}
Now we have just to apply the above Lemma to the first order
formula. Let us write $\la^{(k)}$ for the canonical
representation on $\otimes^k \goth g_{-1}^*\otimes V_\la$.
$$\align
((\nabla^\om)^2\o p^* &- p^*(\nabla^{\ga})^2)(u)(X,Y)=\\
&=\nabla^\om_Y(\la([\_,\ta(u)])\o(p^*s)(u))(X) +
(\la^{(1)}([Y,\ta(u)])(p^*\nabla^{\ga}s)(u))(X)\\
&=\la([X, \Ga(u).Y])(p^*s(u))+\\
&\newquad\la^{(1)}([Y,\ta(u)])(\la([\_,\ta(u)])(p^*s)(u))(X) - \\
&\newquad\tfrac12\la([X,[\ta(u),[\ta(u),Y]]])(p^*s)(u)+\\
&\newquad\la([X,\ta(u)])(p^*(\nabla^{\ga}_Ys))(u) +\\
&\newquad (\la^{(1)}([Y,\ta(u)])(p^*(\nabla^\ga s))(u))(X)
\endalign$$
Altogether we have got

\proclaim{\nmb.{4.6}. Proposition}
For each admissible Cartan connection $\om$, $B_{0}$-equivariant
section $\si$ and for each $B_0$-equivariant function $s\:P_0\to V_{\la}$
$$\align
((\nabla^\om)^2\o p^* &- p^*\o (\nabla^{\ga})^2)s(u)(X,Y)=
\la([X,\Ga(u).Y])(p^*s(u))+\\
&\newquad \la^{(1)}([Y,\ta(u)])(\la([\_,\ta(u)])(p^*s)(u))(X) -\\
&\newquad\tfrac12\la([X,[\ta(u),[\ta(u),Y]]])(p^*s)(u)+\\
&\newquad \la([X,\ta(u)])(p^*(\nabla^{\ga}_Ys))(u)+\\
&\newquad (\la^{(1)}([Y,\ta(u)])(p^*(\nabla^\ga )s)(u))(X)
\endalign$$
holds for all $u\in P$. In particular, vanishing of $\ta$ yields
$$
((\nabla^\om)^2\o p^* - p^*\o (\nabla^{\ga})^2)s(u)(X,Y)=
\la([X,\Ga(u).Y])(p^*s(u)).
$$
\endproclaim

\subhead \nmb.{4.7}. The third order \endsubhead
Exactly in the same way, we use the second order formula to
compute the next one. Let us write briefly
$\text{ad}_{\ta(u)}^2X:= [\ta(u),[\ta(u),X]]$. Furthermore, we shall
write the arguments $X_i$ on the places where they have to be
evaluated, the order of the evaluation is clear from the
context. In fact, whenever $\la^{(k)}$ appears, the arguments
$X_1,\dots, X_k$ are evaluated after this action. We obtain
$$\align
((\nabla^{\om})^3\o p^* &-
p^*(\nabla^\ga)^3)s(u)(X_1,X_2,X_3) =  \\
&=\nabla^{\om}_{X_3}(\text{2nd order difference}) +
\ze_{[\ta(u),X_3]}(\si(p(u))).(p^*(\nabla^\ga)^2s)(X_1,X_2)\\
&=\la^{(2)}([X_3,\ta(u)])\la([X_1,\Ga(u).X_2])(p^*s)(u)+\\
&\newquad\la([X_1,(\nabla^{\ga}_{X_3}\Ga)(u).X_2])(p^*s(u)) +\\
&\newquad \la([X_1,\Ga(u).X_2])((p^*\nabla^{\ga}_{X_3}s)(u))+\\
\noalign{\vskip-2mm\hbox{\hskip4cm\vrule height.2pt
width1.5cm}\vskip-1mm}
&\newquad\la^{(2)}([X_3,\ta(u)])\la^{(1)}([X_2,\ta(u)])
\la([X_1,\ta(u)])(p^*s)(u)-\\
&\newquad\tfrac12\la^{(1)}([X_2,\ta(u)])\la([X_1,\text{ad}_{\ta(u)}^2X_3])
(p^*s)(u)-\\
&\newquad\tfrac12\la^{(1)}([X_2,\text{ad}_{\ta(u)}^2X_3])\la([X_1,\ta(u)])(p^*s)
(u)+\\
&\newquad\la^{(1)}([X_2,\ta(u)])\la([X_1,\ta(u)])
(p^*\nabla^{\ga}_{X_3}s)(u)+\\
&\newquad\la^{(1)}([X_2,\Ga(u).X_3])\o\la([X_1,\ta(u)])(p^*s)(u)+\\
&\newquad\la^{(1)}([X_2,\ta(u)])\o\la([X_1,\Ga(u).X_3])(p^*s)(u)-\\
\noalign{\vskip-2mm\hbox{\hskip4cm\vrule height.2pt
width1.5cm}\vskip-1mm}
&\newquad\tfrac12\la^{(2)}(X_3,\ta(u))\la([X_1,\text{ad}_{\ta(u)}^2X_2])
(p^*s)(u)+\\
&\newquad\tfrac14\la([X_1,[\text{ad}_{\ta(u)}^2X_3,\text{ad}_{\ta(u)}X_2]])
(p^*s)(u)+\\
&\newquad\tfrac14\la([X_1,[\ta(u),[\text{ad}_{\ta(u)}^2X_3,X_2]]])
(p^*s)(u)-\\
&\newquad\tfrac12\la([X_1,\text{ad}_{\ta(u)}^2X_2])
(p^*\nabla^\ga_{X_3}s)(u)-\\
&\newquad\tfrac12\la([X_1,[\Ga(u).X_3,[\ta(u),X_2]]])(p^*s)(u)-\\
&\newquad\tfrac12\la([X_1,[\ta(u),[\Ga(u).X_3,X_2]]])(p^*s)(u)+\\
\noalign{\vskip-2mm\hbox{\hskip4cm\vrule height.2pt
width1.5cm}\vskip-1mm}
&\newquad\la^{(2)}([X_3,\ta(u)])\la([X_1,\ta(u)])(p^*\nabla^\ga_{X_2}s)(u)-\\
&\newquad\tfrac12\la([X_1,\text{ad}_{\ta(u)}^2X_3])
(p^*(\nabla^\ga_{X_2}s))(u)+\\
&\newquad\la([X_1,\ta(u)])p^*(\nabla^\ga_{X_3}\nabla^\ga_{X_2}s)(u)+\\
&\newquad\la([X_1,\Ga(u).X_3])(p^*(\nabla^{\ga}_{X_2}s))(u)+\\
\noalign{\vskip-2mm\hbox{\hskip4cm\vrule height.2pt
width1.5cm}\vskip-1mm}
&\newquad\la^{(2)}([X_3,\ta(u)])\la^{(1)}([X_2,\ta(u)])
(p^*\nabla^\ga_{X_1}s)(u)-\\
&\newquad\tfrac12\la^{(1)}([X_2,\text{ad}_{\ta(u)}^2X_3])
(p^*\nabla^\ga_{X_1} s)(u)+\\
&\newquad\la^{(1)}([X_2,\ta(u)])(p^*(\nabla^\ga_{X_3}\nabla^\ga_{X_1}
s)(u))+\\
&\newquad\la^{(1)}([X_2,\Ga(u).X_3])(p^*(\nabla^\ga_{X_1}
s)(u))+\\
\noalign{\vskip-2mm\hbox{\hskip4cm\vrule height.2pt
width1.5cm}\vskip-1mm}
&\newquad\la^{(2)}([X_3,\ta(u)])(p^*(\nabla^\ga)^2s(u))(X_1,X_2)
\endalign$$
where the horizontal rules indicate the relation to
the individual terms in the second
order difference. Collecting the terms without $\ta$ we obtain
the universal formula for the third order correction terms.

\proclaim{\nmb.{4.8}. Proposition}
For each admissible Cartan connection $\om$, $B_{0}$-equivariant
section $\si$ for each function $s\:P_0\to V_{\la}$ and for all
$u\in \si(P_0)$ we have
$$\align
((\nabla^\om)^3\o p^* - p^*\o (\nabla^{\ga})^3)s(u)(X,Y,Z)
&=\la([X,(\nabla^{\ga}_Z\Ga)(p(u)).Y])(s(p(u))) +\\
&\newquad \la([X,\Ga(p(u)).Y])((\nabla^{\ga}_Zs)(p(u)))+\\
&\newquad\la([X,\Ga(p(u)).Z])((\nabla^{\ga}_Ys)(p(u)))+\\
&\newquad(\la^{(1)}([Y,\Ga(p(u)).Z])((\nabla^\ga s(p(u))))(X).
\endalign$$
\endproclaim

\subhead \nmb.{4.9}. Higher orders\endsubhead
We have seen that the computation of the full formulae goes
quickly out of hands, but it is algorithmic enough to be a good
task for computers.

\proclaim{Algorithm} The difference $F^k s := (\nabla^\om)^k (p^*s)
- p^*((\nabla^\ga)^k s)$ is given by the recursive formula
$$\align
F^0s(u) &= 0\\
F^ks(u)(\row X1k) &= \la^{(k-1)}([X_k,\ta(u)])(F^{k-1}s(u))(\row
X1{k-1}) + \\
&\newquad S_{\ta}(F^{k-1}s(u))(\row X1{k-1}) +\\
&\newquad S_{\nabla}(F^{k-1}s(u))(\row
X1{k-1})+\\
&\newquad S_{\Ga}(F^{k-1}s(u))(\row X1{k-1}) +\\
&\newquad\la^{(k-1)}([X_k,\ta(u)])(p^*((\nabla^\ga)^{k-1}s)(u))(\row
X1{k-1}).
\endalign$$
This expression expands into a sum of terms of the form
$$
a\la^{(t_1)}(\be_1)\dots\la^{(t_i)}(\be_i)p^*(\nabla^\ga)^js
$$
where $a$ is a scalar coefficient,
the $\be_\ell$ are iterated brackets involving
some arguments $X_\ell$, the iterated invariant
differentials $(\nabla^{\ga})^r\Ga$ evaluated
on some arguments $X_\ell$, and $\ta$. Exactly the first $t_j$ arguments
$\row X1{t_j}$ are evaluated after the action of
$\la^{(t_j)}(\be_j)$, the other ones appearing on the right are
evaluated  before. The individual transformations
in $S_\ta$, $S_\nabla$ and $S_\Ga$ act as follows.
\roster
\item The action of $S_\ta$ replaces each summand
$a\la^{(t_1)}(\be_1)\dots\la^{(t_i)}(\be_i)p^*(\nabla^\ga)^js$ by
a sum with just one term for each occurrence of $\ta$ where this
$\ta$ is replaced by $[\ta,[\ta,X_k]]$ and the coefficient $a$ is
multiplied by $-1/2$.
\item $S_\nabla$  replaces each summand in $F^{k-1}$ by a sum
with just one term for each occurrence of $\Ga$ and its
differentials, where these arguments are
replaced by their covariant derivatives $\nabla^\ga_{X_k}$, and
with one additional term where $(\nabla^\ga)^js$ is replaced by
$\nabla^\ga_{X_k}((\nabla^\ga)^js)$.
\item $S_\Ga$ replaces each summand by a sum with just one term
for each occurrence of $\ta$ where this $\ta$ is replaced by
$\Ga(u).X_k$.
\endroster
If we want to compute the
correction terms in order $k$, then during the expansion of
$F^{k-\ell}$ we can omit all terms which involve more then $\ell$
occurrences of $\ta$.
\footnote{
\noindent Some formulae were computed using MAPLEV2.
The number of terms
in low order formulae are\newline
\centerline{
\vbox{\halign{#\hfil\quad&\quad\hfil#\quad&\quad\hfil#\quad&
\quad\hfil#\quad&\quad\hfil#\quad&\quad\hfil#\quad&\quad\hfil#\quad\cr
Order & 1&2&3&4&5&6\cr
Full formula & 1&5&24&134&900&7184\cr
Correction terms & 0&1&4&16&67&328\cr
Linear obstruction terms & 1&2&8&30&153&830\cr
}}}\newline
}
\endproclaim

\demo{Proof}
The algorithm is fully based on Lemma \nmb!{4.4} and the initial
discussion in \nmb!{4.1}.
The last term uses just the equivariancy of the $(k-1)$st covariant
differential. All the other terms correspond exactly to the four
groups of terms in Lemma \nmb!{4.4}. Since we have proved
already in \nmb!{4.4} that an application of this lemma brings
always sums of terms with the required equivariancy properties,
it remains only to verify that the rules deduced in \nmb!{4.4}
yield exactly our formulae.

The first two terms are just in the form derived in \nmb!{4.4}.
The third one is obtained by the differentiation of the induced
mapping $\tilde f$ defined on $P_0$. But this
means that we have to differentiate it like a matrix valued
function, i.e\. we can first evaluate in $\ta$'s and then
differentiate them as constants. Since the whole expression is
multilinear in the arguments involving $\Ga$, the final form of the
transformation follows from the chain rule. The fourth term is
also precisely that one from \nmb!{4.4}.
\qed\enddemo

\subheading{\nmb.{4.10}}
Let us now consider the sections $s\in C^\infty (P, V_\la)^B$ for an
irreducible $B$-representation $\la$ as before,
another irreducible $B$-representation space $V_\mu$, and a linear zero
order operator $\Ph\in\operatorname{Hom}(\otimes^k\goth g_{-1}^*\otimes
V_\la, V_\mu)^{B_0}$. Our
formula for the iterated invariant differential yields
$$
\Ph\o(\nabla^\om)^ks=\Ph\o(\nabla^\ga)^ks + D_0(\ga, \Ga)s +
D_1(\ga,\Ga,\ta)s +\dots +D_k(\ga,\Ga,\ta)s
$$
where $D_j$ collects just those terms
which involve precisely $j$ occurrences of $\ta$.
We call $D_0$ the \idx{\it correction term} while $D_j$,
$j>1$, are called the \idx{\it obstruction terms} of degree $j$ (they are
$j$-linear in $\ta$).
Let us underline, that the correction terms and the obstruction
terms are built by the universal recursive formula based on
\nmb!{4.4}, by means of the same linear mapping $\Ph$.
Their values depend on the initial choice of the equivariant
section $\si\:P_0\to P$, however they turn out to be universal
polynomial expressions in $\nabla^\ga$, $\Ga$ and $\ta$ (but $\ta$
itself depends on the chosen $\si$). Of course, the
composition $D=\Ph\o(\nabla^\om)^k$ is a differential operator transforming
$C^\infty(P,V_\la)^B$ into $C^\infty(P,V_\mu)^B$ if and only if all obstruction
terms vanish independently of the choice of $\si$.

\proclaim{Lemma} The obstruction terms $D_1,\dots,D_k$
vanish for all choices of equivariant sections $\si\:P_0\to P$ if
and only if the first degree obstruction term $D_1$ vanishes for all
choices of $\si$.
\endproclaim
\demo{Proof} Let us
consider a $B_0$--homomorphism $\Ph\:\otimes^k\goth g_{-1}^*\otimes
V_\la\to V_\mu$. The obstruction terms vanish for all choices of $\si$
if and only if $\Ph\o(\nabla^\om)^k(p^*s)$ is $B_1$-invariant for
all $s\in C^\infty(P_0,V_\la)^{B_0}$. This is
equivalent to the vanishing of the derivative
$\ze_{Z}(u).(\Ph\o(\nabla^\om)^k)(p^*s)$ for all $Z\in \goth g_1$ and
$u\in P$. Let us fix $u_0\in P$ and the section $\si\:P_0\to P$
with $u_0\in\si(P_0)$, set $\ga=
\si^*\th_0$, and let $\Ga$ be the (unique) deformation tensor
corresponding to $\tilde\ga$ and $\om$.
Then each of the obstruction terms is expressed as
$$
D_j(\ga,\Ga,\ta)s(u)= f_j(u) = \tilde
f_j(p(u))(\ta(u),\dots,\ta(u))
$$
where $f_j(p(u))\in S^j\goth g_{1}^*\otimes V$ is a homogeneous
polynomial mapping with values in $V_\mu$. Of course, $\tilde f_j\o
p$ are constant in the $\goth g_1$ directions and according to our
choices $\ta(u_0)=0$. Now we can compute
$$\align
\ze_{Z}(u_0).&(\Ph\o(\nabla^\om)^k(p^*s))=\\
&=\ze_Z(u_0).p^*(\Ph\o(\nabla^\ga)^ks + \tilde f_0) + \sum_{j=1}^k
\ze_Z(u_0).((\tilde f_j\o p)(\ta,\dots,\ta))\\
&=(\tilde f_1\o p)(Z)
\endalign$$
where the last equality follows from Lemma \nmb!{4.3}. Thus if
the first degree obstruction terms vanishes for all choices of
$\si$, then $\Ph\o(\nabla^\om)^k(p^*s)$ is $B_1$-invariant as
required.
\qed\enddemo

\subhead\nmb.{4.11}. Remark\endsubhead The above calculus for admissible
connections gives a unified way how to compute the variation of an
expression given in terms of covariant derivatives with respect
to $\ga_0$ and its
curvature tensor, caused by the replacement of $\ga_0$ by another connection
$\ga$ from the distinguished class. Let $\si_0$ and $\si$ be the
$B_0$--equivariant sections corresponding to the connections $\ga_0$ and
$\ga$. Then there is the one form $\ups\in C^\infty(P_0,\goth g_1)^{B_0}$
defined by $\si(u)=\si_0(u).\operatorname{exp}\ups(u)$, see \nmb!{3.6}. Now,
we can use the calculus for the admissible Cartan connections to compare the
induced connections $\tilde \ga_0$ and $\tilde \ga$ and it turns out that
the above expansions in terms of the covariant derivatives and $\ta$'s yield
exactly the variations of the covariant parts. We shall only comment on
these topics here, the details are worked out in \cite{Slov\'ak, 95}.

The relation between the covariant derivatives is
$$
\nabla ^{\ga} s(u)(X) = \nabla^{\ga_0}s(u)(X) + \la([X,\ups(u)])\o s(u),
\quad X\in\goth g_{-1},\ u\in P_0
$$
while the change of the deformation tensors $\Ga_0$ and $\Ga$
transforming $\tilde\ga_0$ and $\tilde\ga$ into another fixed admissible Cartan
connection $\om$ (e.g\. the canonical one) is
$$
\Ga(u)(X)= \Ga_0(u)(X) - \nabla^{\ga_0}\ups(u)(X) -
\tfrac12[\ups(u),[\ups(u),X]], \quad X\in \goth g_{-1},\ u\in P_0
.$$
For the proof see \cite{Slov\'ak, 95, Theorem 1}.

The covariant part of the expansion of $(\nabla^\om)^ks$ in terms of the
connection $\ga_0$ is $(\nabla^{\ga_0})^k s + D_0(\ga_0,\Ga_0)s$.
In terms of the other connection $\ga$, the evaluation on the section
$\si_0$ yields
$$\align
(\nabla^\om)^k &p^{*}s(\si_0(u)) = (\nabla^{\ga_0})^ks(u) + D_0(\Ga_0,\ga_0)s(u) + 0\\
&= (\nabla^{\ga})^ks(u) + D_0(\Ga,\ga)s(u) + D_1(\Ga,\ga,\ups)s(u) +\dots +
D_k(\Ga,\ga,\ups)s(u)
.\endalign$$
Thus, the difference of the covariant parts (i.e\. the variation of this
expression under the change of the underlying connection) is exactly  the
sum of the obstruction terms with $\ups$ substituted for $\ta$. The striking
consequence of this observation is that the covariant parts of the
expansions of the iterated differentials do not involve any derivatives of
the $\ups$'s in their variations under the change of the connection.

Moreover, we can apply our formulae to any linear combination
$D=\sum_{\ell=1}^k A_\ell\o(\nabla^\om)^\ell$ where the zero order operators
$A_\ell$ may be allowed to depend on the curvature of the canonical Cartan
connection $\om$ and its iterated invariant derivatives. Such an expression
defines a natural differential operator if and only if all the obstruction
terms vanish. A more detailed discussion based on our recurrence procedure
and \nmb!{3.10}, \nmb!{3.8}, \nmb!{3.6} shows that we can find all
differential operators built of the covariant derivatives and the curvatures
of the underlying connections $\ga$ and independent on the particular choice
of $\ga$ in this way, see \cite{Slov\'ak, 95, Theorem 2}.

Thus, our procedure extends the methods due to W\"unsch
and G\"unther (developed originally for conformal Riemannian manifolds of
dimensions $m\ge 4$) to all AHS structures. More discussion on various links
to classical methods can be found in \cite{Slov\'ak, 95}.

\head\nmb0{5}. Invariant jets and natural operators \endhead

In this section we discuss the concepts of natural bundles and natural
operators on manifolds equipped with $B$--structures. We show how to
interpret invariant derivatives with respect to Cartan connections as
sections of bundles, and how to naturally construct operators from them.
Since there are canonical Cartan connections on the AHS structures,
this will lead to natural operators.

\subheading{\nmb.{5.1}. Natural bundles and operators}
We shall write $\Mf_m(G)$ for the category of
$m$-dimensional manifolds with almost Hermitian symmetric structures
corresponding to the Lie group $G$ with Lie algebra $\goth g=\goth
g_{-1}\oplus \goth g_0\oplus \goth g_1$ as defined in \nmb!{3.4}.
The morphisms in the category $\Mf_m(G)$
are just principal bundle homomorphisms which cover locally invertible
smooth maps between the bases and preserve the soldering forms.

For each representation $\la\:B\to GL(V_\la)$
and each object $(P,M,\th)\in \Mf_m(G)$ there is the
associated vector bundle $E_\la M$ to the principal bundle $P\to M$
defining the $B$--structure on $M$. This
construction is functorial, we obtain the so called \idx{\it
natural vector bundle\/} $E_{\la}$ on $\Mf_m(G)$. Classically, one
is mainly interested in representations of the first order part $B_0$,
which are trivially extended to the whole $B$. We will devote special
attention to this case, too.

A {\it natural operator} $D\: E_\la \to E_\mu$ between two
natural vector bundles is a system of
operators $D_M\: C^\infty(E_\la M) \to C^\infty(E_\mu M)$
such that for all morphisms $f$ covering a smooth map
$\underline f\: M\to N$ and sections $s_1$,
$s_2\in C^\infty(E_\la M)$, the right-hand square commutes whenever the
left-hand one does
$$\CD
E_\la M
@<s_1<< 
M
@>Ds_1>>
E_\mu M \\
@V{E_\la f}VV
@V{\underline f}VV
@V{E_\mu}VV
\\
E_\la N
@<s_2<<
N
@>Ds_2>>
E_\mu N
\endCD$$

Notice, that the latter definition implies the locality of all
operators $D_M$. A general approach to natural bundles and
operators is developed in \cite{Kol\'a\v r, Michor, Slov\'ak,
93}.

These general definitions of natural bundles and natural operators work well
for each category of manifolds with structures of a fixed type, however, in
our cases the naturality requirements are very weak. The reason is, that
there are nearly no morphisms on general manifolds with AHS structures.
Thus, a stronger restriction of the class of operators under study is
specified by most authors. Mostly one is interested in operators built from
the distinguished linear connections and their curvatures by means of the
covariant derivatives which are independent of any particular choice. Such
operators are usually called \idx{\it invariant} and obviously they are
natural in the above sense.

\subheading{\nmb.{5.2}. The homogeneous case} There is the
subcategory $\Mf_m^{\text{flat}}(G)\subset \Mf_m(G)$ of
spaces locally isomorphic to the homogeneous space $M=G/B$. We can
apply the same definition of the natural operators to this
subcategory. Due to the homogeneity of the objects, each natural
operator on $\Mf_m^{\text{flat}}(G)$ is completely determined
by $D_{G/B}$ and the latter is in turn determined by its action on
germs of sections in one point of $G/B$. The action of the automorphisms
of $G/B$ on the corresponding structure bundle $G\to G/B$ is given by
the left multiplication by the individual elements of $G$. The
sections of the bundles $E_\la(G/B)$ are identified with
$B$-equivariant $V_\la$-valued functions on $G$ and the induced
action of the automorphisms is just the composition of these
functions with the left multiplications by the inverse element.
Thus the operators $D_{G/B}$ with
the invariance properties of our natural operators are exactly the
so called  translational invariant operators, cf\. \cite{Baston,
90}.

This observation suggests another problem on the invariant operators:
What are the invariant operators whose restrictions to the locally flat
spaces coincide with a given natural operator on $\Cal M
f^{\operatorname{flat}}_m(G)$?

The invariant derivatives are manifestly natural operations depending on the
Cartan connection, but they do not map sections of bundles to sections of
bundles. However, we shall build a modification of the standard jet
prolongation and this will lead to operators depending naturally on a Cartan
connection which will play the role of the universal invariant $k$th order
operators.

\subhead \nmb.{5.3}. The jet prolongation of a
representation\endsubhead
As noted above, we would like to view the invariant derivatives of
a section of a natural bundle again as sections of natural bundles. For
the individual derivatives this is impossible, but we can define some
sort of jets. The general idea is to use the invariant differential to
identify the standard first jet prolongation of a natural bundle $E_\la$
with the associated bundle induced by an appropriate $B$-representation.
In fact, this can be done in the general setting, where
$\goth g=\goth g_{-1}\oplus \goth b$ is any Lie algebra, which
linearly splits into the direct sum of an abelian subalgebra $\goth
g_{-1}$ and a subalgebra $\goth b$. So we return to the setting of
chapter \nmb!{2} for the next three subsections.

Assume that we have given a principal $B$--bundle $P\to M$ with a
Cartan connection $\om=\om_{-1}\oplus \om_{\goth b}\in\Om^1(P,\goth
g)$. Moreover assume that
$\la:B\to \text{GL}(V_\la)$ is a representation of $B$ and $s\in
C^\infty(P,V_\la)$ is a smooth map, and consider the smooth map
$(s,\nabla^\om s):P\to V_\la\oplus (\goth g_{-1}^*\otimes V_\la)$. Then
for each $Z\in \goth b$, we obtain
$$\align
\ze_{Z}.(s,\nabla^{\om}_Xs) &=
(\L_{\om^{-1}(Z)}s,\L_{\om^{-1}(Z)}\o\L_{\om^{-1}(X)}s) \\
&=(\ze_Z.s,\nabla^\om_X(\ze_Z.s) + \L_{\om^{-1}([Z,X])}.s),
\endalign$$
where we have essentially used the horizontality of the curvature
of any Cartan connection. Assume now that $s\in
C^\infty(P,V_\la)^B$
is equivariant. Then $\ze_Z.s=-\la(Z)\o s$ and we get
$$
-\ze_{Z}.(s,\nabla^{\om}_Xs) = (\la(Z)\o
s,\la(Z)\o(\nabla^\om_Xs) -\nabla^\om_{[Z,X]_{-1}}s +
\la([Z,X]_{\goth b})\o s),
$$
where we have split $[Z,X]=[Z,X]_{-1}+[Z,X]_{\goth b}$ according to
the decomposition of $\goth g$.

Thus we define the space
$\J^1(V_\la) := V_\la\oplus (\goth g_{-1}^*\otimes V_\la)$
and the mapping $\tilde\la\:\goth b\x\J^1(V_\la)\to \J^1(V_\la)$
by the formula
$$
\tilde\la(Z)(v,\ph) =(\la(Z)(v),
\la(Z)\o\ph -\ph\o\ad_{-1}(Z) + \la(\ad_{\goth b}(Z)(\_))(v))
$$
where $\ad_{-1}(Z):\goth g_{-1}\to \goth g_{-1}$ is the map $X\mapsto
[Z,X]_{-1}$ and $\la(\ad_{\goth b}(Z)(\_))(v)\:\goth g_{-1}\to V_\la$
is defined by $\la(\ad_{\goth b}(Z)(\_))(v)(X)=\la([Z,X]_{\goth b})(v)$.

\proclaim{Lemma} The mapping $\tilde \la$ is an action of $\goth
b$ on $\J^1(V_\la)$. For each $\goth b$--equivariant element
$s\in C^\infty(P,V_\la)$,  the mapping
$(s,\nabla^\om s)\: P\to \J^1(V_{\la})$ is
$\goth b$-equivariant with respect to this action.
\endproclaim
\demo{Proof} For $Z$, $W\in \goth b$ and $(v,\ph)\in\J^1(V_{\la})$ we
compute:
$$\align
\tilde\la(W)&\tilde\la(Z)(v,\ph)=(\la(W)\la(Z)(v),\la(W)\o\la(Z)\o\ph-
\la(W)\o\ph\o\ad_{-1}(Z)+\\
&+\la(W)\o\la(\ad_{\goth b}(Z)(\_))(v)-\la(Z)\o\ph\o\ad_{-1}(W)+
\ph\o\ad_{-1}(Z)\o\ad_{-1}(W)-\\
&-\la(\ad_{\goth b}(Z)(\_))(v)\o\ad_{-1}(W)+
\la(\ad_{\goth b}(W)(\_))(\la(Z)(v)))
.\endalign$$
Now when forming the commutator of $\tilde\la(W)$ and $\tilde\la(Z)$
the second and fourth term in the second component do not contribute,
so we get
$$
(\tilde\la(W)\tilde\la(Z)-\tilde\la(Z)\tilde\la(W))(v,\ph)=(\la([W,Z])(v),
\la([W,Z])\o\ph+\Ph),
$$
where $\Ph$ is the linear map defined by
$$\align
\Ph(X)&=\la(W)\la([Z,X]_{\goth b})(v)-\la(Z)\la([W,X]_{\goth b})(v)+
\ph([Z,[W,X]_{-1}]_{-1})-\\
&-\ph([W,[Z,X]_{-1}]_{-1})-\la([Z,[W,X]_{-1}]_{\goth b})(v)+
\la([W,[Z,X]_{-1}]_{\goth b})(v)+\\
&+\la([W,X]_{\goth b})\la(Z)(v)-\la([Z,X]_{\goth b})\la(W)(v)
.\endalign$$
Now using the Jacobi identity and the fact that $\goth b$ is a subalgebra
while $\goth g_{-1}$ is abelian, one immediately verifies that
$$\gather
[[W,Z],X]_{-1}=[W,[Z,X]_{-1}]_{-1}-[Z,[W,X]_{-1}]_{-1}\\
{\aligned
[[W,Z],X]_{\goth b}&=[W,[Z,X]_{\goth b}]-[Z,[W,X]_{\goth b}]+\\
&+[W,[Z,X]_{-1}]_{\goth b}-[Z,[W,X]_{-1}]_{\goth b}.
\endaligned}
\endgather$$
Using this one immediately sees that
$$
\Ph(X)=-\ph([[W,Z],X]_{-1})+\la([[W,Z],X]_{\goth b})(v).
$$
The rest of the lemma is a consequence of our definition.
\qed\enddemo

Thus we can consider the mapping $s\mapsto(s,\nabla^\om s)$ as a
section of the associated bundle to $P$ induced by the $\goth b$-module
$\Cal J^1V_\la$.

In fact, the action $\tilde \la$ coincides with the canonical
action of $\goth b$ on the standard fiber of the usual first jets of
sections of $E_\la(G/B)$ (which also could be used as a geometrical argument
for the proof of the above Lemma). Thus our construction can be understood
as identifications of the standard first jet prolongations $J^1E_\la M$ with
the associated bundles $P\x_B \Cal J^1V_\la$, determined by the Cartan
connection $\om$. The section $(s,\nabla^\om s)$ of this bundle is then
called the \idx{\it invariant one--jet} of the section $s$ of $E_\la M$.

\subheading{\nmb.{5.4}. The second jet prolongation of a representation}
Observe that it is easy to extend $\Cal J^1(\_)$ to a
functor on the category of $\frak b$--representations. In fact for a
homomorphism $f\:V\to W$ of $\frak b$--modules we just define
$\Cal J^1(f)\:\Cal J^1(V)\to\Cal J^1(W)$ by
$\Cal J^1(f)(v,\ph):=(f(v),f\o \ph)$. One easily computes directly
that $\Cal J^1(f)$ is a module homomorphism.

Second, it is clear that by projecting onto the first component we get
a module homomorphism $\Cal J^1(V)\to V$, and actually these
homomorphisms constitute a natural transformation between
$\Cal J^1(\_)$ and the identity functor.

Next, let us consider $\Cal J^1(\Cal J^1(V))$. There are two natural
homomorphism from this space to $\Cal J^1(V)$: First, we have the
above mentioned natural projection, and second, there is the first jet
prolongation of the projection $\Cal J^1(V)\to V$, and we define
$\Cal J^2(V)$ to be the submodule of $\Cal J^1(\Cal J^1(V))$ on which
these two homomorphisms coincide. The underlying vector space of
$\Cal J^1(\Cal J^1(V))$ is just
$(V\oplus (\frak g_{-1}^*\otimes V))\oplus\frak
g_{-1}^*\otimes(V\oplus (\frak g_{-1}^*\otimes V))\cong V\oplus
(\frak g_{-1}^*\otimes V)\oplus(\frak g_{-1}^*\otimes V)
\oplus(\frak g_{-1}^*\otimes\frak g_{-1}^*\otimes V)$, and under this
identification $\Cal J^2(V)$ is just the submodule of those elements
where the two middle components are equal. One immediately verifies
that $\Cal J^2(\_)$ is again a functor and that projecting out the
first two components gives a natural transformation to
$\Cal J^1(\_)$.

\subheading{\nmb.{5.5}. Higher jet prolongations}
We can iterate the above procedure as follows: Suppose we have
already constructed functors $\Cal J^i(\_)$ for $i\leq k$ such that
$\Cal J^i(V)$ is a submodule in $\Cal J^1(\Cal J^{i-1}(V))$ and such
that for each $i$ there is a natural transformation
$\Cal J^i(\_)\to\Cal J^{i-1}(\_)$ induced by the projection
$\Cal J^1(\Cal J^{i-1}(V))\to \Cal J^{i-1}(V)$, i.e\. by the natural
transformation from $\Cal J^1(\_)$ to the identity.

Then consider $\Cal J^1(\Cal J^k(V))$ for some module $V$. We
have two natural homomorphisms
$\Cal J^1(\Cal J^k(V))\to \Cal J^1(\Cal J^{k-1}(V))$, namely the
natural projection $\Cal J^1(\Cal J^k(V))\to \Cal J^k(V)$ followed by
the inclusion of the latter space into $\Cal J^1(\Cal J^{k-1}(V))$
and the first jet prolongation of the natural map
$\Cal J^k(V)\to \Cal J^{k-1}(V)$, and we define $\Cal J^{k+1}(V)$ to
be the submodule where these two module homomorphisms coincide.
Moreover for a module homomorphism $f\:V\to W$ we define
$\Cal J^{k+1}(f)$ as the homomorphism induced by $\Cal J^1(\Cal J^k(f))$.
Finally, from the obvious projection $\Cal J^{k+1}(V)\to \Cal J^k(V)$
we clearly get a natural transformation
$\Cal J^{k+1}(\_)\to \Cal J^k(\_)$.

Also by induction, it is easy to see that as a vector space we
always have
$\Cal J^k(V)\cong \oplus_{i=0}^k(\otimes^i\frak g_{-1}^{*}\otimes V)$.
Moreover starting from lemma \nmb!{5.3} it is again clear by
induction that for any $B$--equivariant function $s:P\to V$ and any
Cartan connection $\om$ on $P$, the mapping
$$j^k_\om s:=(s,\nabla^\om s,\dots ,(\nabla^\om)^k s):P\to \Cal J^k(V)$$
is equivariant, too. This map is called the {\it invariant $k$--jet of
$s$ with respect to $\om$}.

\subheading{\nmb.{5.6}. The AHS case}
Since $\goth g=\goth g_{-1}\oplus \goth g_0\oplus\goth g_1$, there are a few
simplifications in the construction of the jet prolongations. First of all
for $A\in \goth g_0$ we clearly have $\ad_{\goth b}(A)=0$, while for $Z\in
\goth g_1$ we have $\ad_{-1}(Z)=0$, so the action on the first jet
prolongation becomes easier for each case. In particular, we see that in
fact the action of $\goth g_0$ is just the tensorial one. Thus, the
isomorphism $\Cal J^k(V)\cong \oplus_{i=0}^k(\otimes^i\frak
g_{-1}^{*}\otimes V)$ is not only an isomorphism of vector spaces but also
of $\goth g_0$--modules.

Moreover, in this case we have the additional information that the
group $B$ is the semidirect product of the contractible subgroup $B_1$
and the subgroup $B_0$. If we start with a $B$--representation and
form the first jet prolongation of the corresponding $\goth
b$--representation, then this will always integrate to a
$B$--representation. This is due to the fact that the restriction to
$\goth g_1$ integrates by contractibility of $B_1$, while the action
of $\goth g_0$ is the tensor product of the original action with the
adjoint action, and both of these integrate.

Therefore in this case, for each representation $\la:B\to GL(V_\la)$
we obtain the jet prolongations $\Cal J^k(\la):B\to GL(\Cal
J^k(V_\la))$. This in turn implies that for each natural bundle
$E_\la$ on $\Mf_m(G)$ there is the natural bundle $\Cal J^k(E_\la)$.
By the construction, this bundle coincides with the so
called $k$th semi-holonomic jet prolongation of $E_\la$.

\subheading{\nmb.{5.7}}
There is now a simple procedure how to use the invariant jets with respect
to a Cartan connection for the constructions of
differential operators.  Suppose that
$\la\:B\to GL(V_\la)$, is a representation of $B$, and
suppose that for some $k$ and
another such representation $\mu$ on $V_\mu$, there is a $B$--equivariant
(even nonlinear) mapping $\Ph\:\Cal J^k(V_\la)\to V_\mu$.
Now, for each $P\to M$ with a Cartan connection $\om$, we can define a
$k$--th order differential operator
$D_M:C^\infty(E_\la M)\to C^\infty(E_\mu M)$ by putting
$$D_M(s)(u):=\Ph(j^k_\om s(u))\qquad \forall s\in C^\infty(P,V)^B,
u\in P.$$
The associated bundles $E_\la M$, $E_\mu M$ are functorial in $P\to M$,
and so by the construction the
operators defined in this way intertwine the actions of all morphisms
of the $B$--structures on the sections, which preserve the Cartan
connections. In particular, if there is a canonical Cartan connection,
which is preserved under the action of all morphisms, then the
operator constructed in this way will be natural. More generally, one
can also interpret these operators as natural operators which also
depend on the Cartan connections, but we will not work out this point
of view here.

\subheading{\nmb.{5.8}} Let us discuss in more detail now, how to find the
$\goth b$-module homomorphisms $\Ph\:\Cal J^k V_\la\to V_\mu$
for irreducible representations $\la$ and $\mu$ of $B_0$ on
$V_\la$ and $V_\mu$, viewed as irreducible representations of $\goth b$.
Let us recall that $\Cal J^k(V)= \oplus_{i=0}^k(\otimes^i\frak
g_{-1}^{*}\otimes V)$ as $\goth g_0$-module.

\proclaim{Lemma}
Let $\pi\:\Cal J^k(V_\la)\to
\otimes ^k\frak
g_{-1}^{*}\otimes V_\la$ be the $\frak g_0$--homomorphism
corresponding to the decomposition of $\Cal J^k(V_\la)$ and
let $\Ph\:\Cal J^k V_\la\to V_\mu$ be a $\goth g_0$-module
homomorphism whose restriction to $\otimes^k\frak
g_{-1}^{*}\otimes V\subset \Cal J^kV_\la$ does not vanish.
Then $\Ph$ is a $\goth b$-module homomorphism if and only if it
factors through $\pi$ and $\Ph$ vanishes on the image of $\otimes ^{k-1}
\frak g_{-1}^{*}\otimes V_\la\subset \Cal J^k(V_\la)$
under the action of $\goth b_1$.
\endproclaim
\demo{Proof}
Let $\Bbb I$ be a generator of the center of $\frak g_0$. Then by
Schur's lemma $\Bbb I$ acts by a scalar on every irreducible
representation of $\frak g_0$. Moreover, for the adjoint representation
these scalars are just given by the grading. Now the action of $\frak
g_0$ on each of the components $\pi^k_i(\Cal J^k(V_\la))$ is the
tensorial one, so $\Bbb I$ acts by different scalars on each of them.
Moreover, any $\goth g_0$-module homomorphism $\Ph$ defined on the top part
of $\Cal J^k(V_\la)$ is a $\goth b$-module homomorphism if and only if $\Ph$
is $\goth g_1$-invariant.
\qed
\enddemo

\proclaim{\nmb.{5.9}. Lemma}
The action of an element $Z\in \goth g_1$ on
$Y_1\otimes\dots\otimes Y_{k-1}\otimes v_{k-1}\in \otimes ^{k-1}\frak
g_{-1}^{*}\otimes V_\la\subset \Cal J^kV_\la$ yields
$$
\sum_{i=0}^{k-1}\biggl(\sum_\al Y_1\otimes\dots\otimes Y_i\otimes\et_\al\otimes
\bigl([Z,\xi_\al].(Y_{i+1}\otimes\dots\otimes Y_{k-1}\otimes v_{k-1})
\bigr)\biggr) \in \otimes ^k\frak g_{-1}^{*}\otimes V_\la
$$
where $\et_\al$ and $\xi_\al$ are dual basis of $\goth g_1$ and $\goth g_{-1}$
with respect to the Killing form and the dot means the canonical action of
the element in $\goth g_0$ on the argument.
\endproclaim
\demo{Proof}
The statement follows easily from the definition of $\Cal J^kV_\la$
by induction on the order $k$.
\qed\enddemo

\subhead\nmb.{5.10}. Remark\endsubhead A reformulation of the preceding
Lemma reads: {\sl A $\goth g_{0}$-module homomorphism
$\Ph\: \otimes^{k}\goth g_{-1}^{*}\otimes V\to W$ can be considered as a
$\goth b$-module homomorphisms $\Cal J^k(V_\la)\to V_\mu$
if and only if
$$
\Ph\bigg(\sum_{i=1}^k\la^{(i-1)}([Z,X_i])\ps
(X_1,\dots, X_{i-1},X_{i+1},\dots,X_k) \bigg)=0$$
for all elements $\ps\in \otimes^{k-1}\goth g_{-1}^{*}\otimes V$, all
$X_1,\dots, X_k\in \goth g_{-1}$, and all $Z\in\goth g_1$.
}

This expression can be also found among the obstruction terms in the
expansion of the $k$th iterated invariant differential $(\nabla^\om)^k$.
Indeed, the linear obstruction terms involve in particular the terms with
highest order derivatives of the section, i.e\. those of order $k-1$ in s,
and a simple check shows that they are of the above form. Let us call this
part the \idx{\it algebraical obstruction term}. Now, the above Lemma
implies that if this algebraical obstruction vanishes `algebraically', i.e\.
before substitution of the values of the invariant jets, then all other
obstruction terms vanish as well and we have got a natural operator in this
way.

Once we have a correspondence between the $\goth b$-module homomorphisms of
the jets and the natural operators, we should try to extend the algebraic
methods leading to  the well known  classification of all linear natural
operators on the locally  flat spaces to the general setting. We shall come
back to this point in the fourth part of this series. A more straightforward
generalization of the Verma module technique can be found in the forthcoming
paper \cite{Eastwood, Slov\'ak}.

An important observation is, that not all operators are created in such
an algebraic way, there are also examples of operators where the algebraic
obstruction does vanish only after the substitution of the invariant jets.
The simplest example is the second power of the Laplacian on the four
dimensional conformal Riemannian manifolds. A more detailed discussion on
such cases can be found in \cite{Eastwood, Slov\'ak}.

\head\nmb0{6}. Remarks on applications\endhead

\subheading{\nmb.{6.1}}
Let us indicate now in more detail
how the theory developed so far applies to the study of natural operators.
The simplest possibility is the one discussed in the end of the previous
section:
\roster
\item Starting with irreducible representations $\la$ and $\mu$
of $B$, we consider all the compositions $\Ph\o(\nabla^\om)^k$, where
$\Ph\:\otimes^k\frak g_{-1}^{*}\otimes V_\la\to V_\mu$
is $\frak g_0$--equivariant and linear.
\item Such an expression yields a differential operator on sections if
and only if $\Ph\o(\nabla^\om)^ks$ is $B_1$--invariant for each
$s\in C^\infty(P,V_\la)^{B}$. In view of the expansion of the iterated
differential in terms of the underlying connections, this
is equivalent to the vanishing of the linear obstruction terms after
substitution of the invariant jets. Moreover, the algebraic vanishing of the
algebraic obstruction terms suffices, see \nmb!{4.10} and \nmb!{5.10}.
\item There are the canonical Cartan connections
$\om$ on all manifolds with AHS
structures  and so the differential
operators obtained in (2) with help of $\om$ turn out to be natural.
\item If we choose a linear connection $\ga$ in the distinguished class, then
there is the unique deformation tensor $\Ga$ transforming the induced Cartan
connection $\tilde \ga$ into the canonical one.
Thus the formulae from Section \nmb!{4}
express the natural operators by means of the covariant
derivatives and curvatures of the linear connection $\ga$.
\endroster
The general construction of the canonical connection $\om$ and the
deformation tensors $\Ga$ are postponed to the next part of the series.
However, in order to have some conrete examples,
we present the computations in the conformal Riemannian case below.

\subheading{\nmb.{6.2}. Conformal Riemannian structures} The existence of
the principal bundle $P\to M$ with a canonical Cartan connection is well
known in the conformal case, see \cite{Kobayashi, 72}.
But we prefer to present an explicit construction  here to illustrate the
links of our concepts and formulae to the the classical approach.

Let us start with a manifold $M$ of dimension $m\geq 3$ equiped with a
conformal class
of Riemannian metrics or equivalently with a reduction of the first
order frame bundle $P^1M$ to the group $B_0=CSO(m)\simeq SO(m)\sempr \Bbb
R$. Any metric in the conformal class has its Levi--Civit\`a
connection, which is torsion free. There is a bijective correspondence
between torsion--free connections on $M$ and $GL(m)$--equivariant
sections of $P^2M\to P^1M$, see \cite{Kobayashi, 72, Proposition
7.1}. Our group $B$ can be viewed as a subgroup of $G^2_m$, see Lemma
\nmb!{3.2}. It turns out that the orbit of the images of the
Levi--Civit\`a connections under the group $B$ coincide, and actually
give a reduction of $P^2M$ to the group $B$, and thus a torsion free
$B$--structure on $M$ by \nmb!{3.5}.

If we start with a reduction  $\ph:P_0\to P^1M$  in the sense of
\nmb!{3.2} to the group $B_0=Spin(m)\sempr \Bbb R$ then one can still
construct a subbundle $\tilde P$ of $P^2M$ as above, and
$P:=\ph^*\tilde P$ is again a torsion free $B$--structure with the induced
soldering form. Thus we can include the Spin representations in our
approach.

To obtain a canonical Cartan connection on such $B$--structures we
proceed as follows: The values of the
$\goth g_0$-component $\ka_0$ of the curvature function $\ka$ of
a Cartan connection $\om$ can be viewed as elements
in $\goth g_{-1}^*\otimes\goth g_{-1}^*\otimes \goth
g_{-1}^*\otimes \goth g_{-1}$, cf\. \nmb!{3.1}.(2). There are three
possible evaluations in the target space. The evaluation
over the last two entries is just the trace in $\goth g_{0}$,
the other two possibilities coincide up to a sign since $\ka$ is a two form.

\subheading{Definition} The \idx{\it
trace of the curvature} $\ka_0$ is the composition of $\ka_0$
with the evaluation over the first and the last entry.
A \idx{\it normal Cartan connection} $\om\in\Om^1(P,\goth g)$
is an admissible connection with a trace-free curvature $\ka_0$.

The general obstruction to the
existence of a normal Cartan connection is in certain cohomology
group, we shall not discuss this point here, cf\. \cite{Ochiai,
70}, \cite{Baston, 90}. But we shall use the
formula \nmb!{3.10}.(4) for the deformation of the curvature in
order to compute explicitly the necessary deformation tensor $\Ga$
for a given admissible connection. It turns out
that the result is uniquely determined by the initial data.

\subhead \nmb.{6.3}\endsubhead
We have to use a coordinate notation for the values of $\Ga$
and $\ka_0$ in order to handle the proper evaluations in the
trace.
So let $e_i$ be the standard basis of the vector space
$\goth g_{-1}$, $e^i$ the standard dual basis in $\goth g_{-1}^*$ and
$e^i_j$ the standard basis of $\goth g_{-1}^*\otimes \goth g_{-1}$.
Note that the bases $e_i$ and $e^i$ are in fact dual with respect to
the Killing form, up to a fixed scalar multiple. Then
$\Ga(u)(e_{i}) = \sum_j\Ga_{ji}(u).e^j$,
$\ka_0(u)(e_{i},e_j)= \sum_{k,l} K^k_{lij}(u).e^l_k$. In the
sequel, we shall not always indicate explicitly sums over repeated
indices. If we restrict the manipulations with these symbols to
permutations of indices, contractions and similar invariant
tensorial operations, our computations will be manifestly
independent of the choice of the basis. In particular, the
trace of $\ka_0$ is expressed by the functions $K^i_{lij}$.

The brackets of the generators of $\goth s\goth o(m,\Bbb R)$,
$m>2$, are computed
easily from the block-wise representation in \nmb!{3.3}:
$$
[e_i, e^j] = e^j_i-e^i_j + \de^j_i\Bbb I_m,\newquad [e^i_j,e_k]=\de^i_k e_j
$$
where $\Bbb I_m$ stands for the unit matrix. Now we evaluate the formula
for the deformation $\bar\ka_0-\ka_0 =: \de(K^k_{lij})$
of the curvature caused by a choice of $\Ga$, see \nmb!{3.10}.(4).
$$\align
[\Ga.e_j,e_i]-[\Ga.e_i,e_j]&=
\tsum_p\Ga_{pj}(-e^p_i+e^i_p-\de^i_p\Bbb I_m) - \tsum_p
\Ga_{pi}(-e^p_j+e^j_p-\de^j_p\Bbb I_m)\\
&=(-\Ga_{kj}\de^l_i + \Ga_{lj}\de^i_k -\Ga_{ij}\de^l_k +
\Ga_{ki}\de^l_j -\Ga_{li}\de^j_k + \Ga_{ji}\de^l_k)e^k_l
\endalign$$
Thus, the deformation of the trace achieved by $\Ga$ is
$$\align
\de (K^{l}_{klj}) &= (m-3)\Ga_{kj} + \Ga_{jk} +
\de^k_j\tsum_i\Ga_{ii}\\
\de(K^k_{kij}) &= m(\Ga_{ji}-\Ga_{ij})\\
\tsum_j\de (K^i_{jij})&= 2(m-1)\tsum_i\Ga_{ii}.
\endalign$$
We need the third `contraction' for technical
reasons.

Now, assume first we have two normal Cartan connections and let $\Ga$
be the corresponding deformation tensor. Since the torsion is zero,
the  Bianchi identity shows
that for any normal Cartan connection, not only the trace defined in
\nmb!{6.2} but also the trace inside $\goth g_0$ vanishes. Thus the resulting
deformation of all three contractions above must be zero.
So in particular, $\sum_i\Ga_{ii}=0$ and the functions $\Ga_{ij}$
are symmetric in $i$, $j$. But then the first equation yields
$0=(m-2)\Ga_{lj}$. Thus if there is a normal Cartan connection,
it is unique.

Let $\ga$ be the Riemannian connection of an arbitrary metric
from the conformal class on $M$. Then it induces an admissible
connection $\tilde\ga$ on $P$, see \nmb!{3.8}. Moreover, the $\goth
g_{0}$-component of the curvature of the induced connection is
just the pullback of the Riemannian curvature to $P$.
Let us try to deform $\tilde \ga$ by means of symmetric functions
$\Ga_{ij}$.

The deformation is expressed above in the form
$\text{Tr}(\ka_0-\bar \ka_0)$, where $\bar \ka_0$ is the `new one'.
Thus we have just to solve the above equations with respect to
$\Ga$ with $\de(K^k_{lij})$ replaced by $-R^k_{lij}$, the
Riemannian curvature. We obtain easily
$$
\Ga_{ij}=\frac{-1}{m-2}\bigl(R_{ij} -
\frac{\de_{ij}}{2(m-1)}R\bigr),\tag1
$$
where $R_{ij}$ and $R$ are the pullbacks of the Ricci tensor and
the scalar curvature to $P$ (expressed in the frame form, i.e\.
as functions on $P$). Let us notice that the above deformation
tensor $\Ga$ is exactly the so called ``rho--tensor'' used
extensively in conformal geometry because of this ``beautiful
transformation properties''.

Altogether we have reproved, even for conformal Spin structures:

\proclaim{\nmb.{6.4}. Theorem} Let $M$ be a connected smooth manifold,
$\text{\rm dim}M\ge 3$, with a conformal structure $P_0\to M$.
Then there is a unique normal Cartan connection $\om$ on $P\to M$
which is expressed by means of any Riemannian connection $\ga$
from the conformal class by the formula $\om = \tilde \ga -
\Ga\o\th_{-1}$ with $\Ga$ defined by \nmb!{6.3}.(1).
\endproclaim

\subheading{\nmb.{6.5}. Operators on locally flat manifolds}
Now we can apply the canonical normal Cartan connections in the
construction from \nmb!{5.7}. In view of the next lemma, this
procedure yields at least all natural operators ``visible'' on the
locally flat manifolds.

Let us fix two representations $\la$ and $\mu$
of $B_0$ and let $E_\la$ and $E_\mu$ be the corresponding natural
bundles on the manifolds with the conformal (Spin) structures.
Further let us consider the locally flat structures $P\to P_0\to M$.
This means, we assume that there are (locally defined) connections in
the distinguished class with vanishing curvature, or equivalently,
$P\to M$ is locally isomorphic to the homogeneous space $G\to G/B$.

\proclaim{Lemma} Suppose that the family of operators
$D_M\:C^\infty(E_\la M) \to C^\infty(E_\mu M)$ is a
natural operator on the category of locally flat conformal (Spin)
structures and let
$\Pi\o(\nabla^\ga)^k$ be its expression in the (locally defined)
flat connection $\ga$ in the distinguished class on $M$.
Then the operator $\tilde D =\Pi\o (\nabla^\om)^k$ defined by
means of the invariant differential with respect to the unique normal
Cartan connection $\om$ on $P\to M$ transforms $B$-equivariant functions
into $B$-equivariant functions and equals to $D_M$.
\endproclaim
\demo{Proof}
Since the operator $D$ is natural,
$D_M\:C^\infty(PM,V_\la)\to C^\infty(PM,V_\mu)$ commutes with the
induced action of the morphisms which is given by the composition
with the inverses. On the other hand, $\tilde D$ commutes with these
actions as well and since the structure in question is locally flat,
the automorphisms of $P\to M$ act transitively.
Thus, if we show that $\tilde D$ coincides on $PM$
with $D_M$ in one point of $PM$, then they must coincide
globally. But if we choose a flat local connection $\ga$
and the corresponding (local) $B_0$-equivariant section
$\si\:P_0\to P$,
then the unique normal Cartan connection $\om$ equals to the
induced admissible Cartan connection $\tilde \ga$, in particular the
corresponding deformation tensor $\Ga$ is zero. Thus, according to
the preceding section, the iterations of the invariant derivative
with respect to $\om$ and the pullbacks of the iterations of the
covariant derivative with respect to $\ga$ coincide on $\si(P_0)$.
In particular, the operator $\tilde D$ transforms sections into
sections.
\qed\enddemo

\subheading{\nmb.{6.6}. Remark}
By virtue of the general theory of natural operators on  Riemannian
manifolds, the naturality assumption in the previous lemma means just
that the operator $D$ is defined by a universal expression in terms of
the underlying Riemannian connections in the conformal class, see
e.g\. \cite{Kol\'a\v r, Michor, Slov\'ak, 93}. Thus our result shows
that the ``conformally invariant operators'' in the usual sense (see
e.g\. \cite{Branson, 82}, \cite{W\"unsch, 86}, \cite{Baston, Eastwood,
90}) are all obtained by our procedure, at least in the conformally
flat case. Moreover, if we allow more general linear combinations of the
iterated invariant differential (involving the iterated invariant
differentials of the Weyl curvature in dimensions $m\ge 4$, or
the invariant differentials of the Cotton-York tensor for $m-3$),
then we can achieve all the invariant operators mentioned above,
cf\. Remark \nmb!{4.11}.

Furthermore, the lemma is not restricted to linear operators, on the
contrary, the same arguments apply if the expression for the operator
$D_M$ is a polynomial in the covariant derivatives.

\subhead \nmb.{6.7}. Examples \endsubhead
To illustrate the use of the general formulae, let us consider now
some special cases. As before, we shall restrict the attention here to the
conformal case.

Consider an irreducible representation $\la\:B_0\to GL(V_{\la})$
and let us write $\la'$ for its restriction to the semisimple
part of $B_{0}$. Each $\la$ is given by $\la'$ and the scalar
action of the center, $\la(\Bbb I_{m})(v)= -w\cdot v$. The
scalar $w$ is called the \idx{\it conformal weight\/} of $\la$.

According to the above discussion, if there is a
$\goth g_{0}$--homomorphism $\Ph:\otimes^k\goth g_{-1}^*\otimes
V_\la\to V_\rh$ onto an irreducible representation $V_{\rh}$ such
that $\Ph\o(\nabla^{\om})^k$ is a natural operator
then the formulae obtained in Section \nmb!{4}
yield its expression by means of a universal formula in
terms of the underlying linear connections. Recall that we
denoted by $\Ga$
the deformation tensor determined by a choice of a metric in
the conformal class.

We shall look first at the second order operators. For each irreducible
representation $V_{\la}$ of $B_0$, the tensor product
$\goth g_{-1}^*\otimes V_{\la}$ decomposes uniquely into irreducible
representations $V_{\rh}$ (i.e\. there are no multiplicities in the
decomposition), see e.g\. \cite{Fegan, 76}. Let us write
$\operatorname{Id}= \sum_{\rh} \pi^{\la\rh}$ for this decomposition.

Let $V_{\la}$ be an irreducible representation of $\goth g_0$ and
let $\pi^{\la\rh\si}$ be a projection of
$\otimes^2\goth g_{-1}^*\otimes
V_\la$ onto an irreducible representation $V_\si$
given by
$
\pi^{\la\rh\si}(Z_1\otimes Z_2\otimes s)=
\pi^{\rh\si}[Z_1\otimes(\pi^{\la\rh}(Z_2\otimes s))].$
Lemma \nmb!{5.9} gives a possibility to prove that
$\pi^{\la\rh\si}\o(\nabla^{\om})^2$ is a natural
operator for certain choices of $\la$, $\rh$ and $\si$,
and Proposition \nmb!{4.6} is saying that  the natural operator
can be written (using the underlying linear connections) as
$\pi^{\la\rh\si}\{(\nabla^{\ga})^2 +\la([X,\Ga .Y])\}s.$ This is a universal
formula valid for any dimension, any representation and any
projection (even for the other structures, not only for the conformal
one).

Choosing a specific representation, the formula can be simplified
further. Let us consider now for simplicity the case of an even
dimension $m=2k$ and let $e^i, i=1,\ldots,m$ be
weights of the representation $\gog_1.$

{\bf 1.}
 Let us discuss a simple example - second order
operators acting on functions (having possibly  a conformal weight).
Hence let $\la'=0$ be the heighest weight of the trivial
representation $V_{\la'}=\Bbb C$ and $w$ its  conformal weight.
The tensor product $\ox^2\gog_{-1}^*\ox V_{\la}$ decomposes into
three irreducible parts, namely $S^2_0(\gog_{-1}^*)$
(symmetric traceless tensors), the trivial representation and
$\La^2(\gog_{-1}^*)$. Let $\pi_1,\pi_2,\pi_3$ denote the
corresponding projections.

We can use now the algebraic conditions discussed in \nmb!{5.9}.
So $\xi_{\al},$ resp\. $\eta_{\al}$ are dual bases in
$\gog_{-1},$ resp. $\gog_1.$ We have to consider elements of the
form
$$\multline
\sum_{\al}\left\{\eta_{\al}\ox [Z,\xi_{\al}]\cdot(Y\ox v)+
Y\ox \eta_{\al}\ox [Z,\xi_{\al}]\cdot v\right\}=\\
\sum_{\al}\left\{\eta_{\al}\ox([[Z,\xi_{\al}],Y]\ox v)+
w\eta_{\al}\ox Y\ox\xi_{\al}(Z)v+
w Y\ox \eta_{\al}\ox\xi_{\al}(Z) v\right\}.
\endmultline
$$
Using
$[[Z,\xi_{\al}],Y]=-\langle Z,Y\rangle\et_{\al}+\xi_{\al}(Z)Y+\xi_{\al}(Y)Z$
and
$
\sum_{\al}\xi_{\al}(Z)\eta_{\al}=Z,
$
we get
$$
(w+1)[Z\ox Y\ox v+Y\ox Z\ox v]-\langle Z,Y\rangle(\sum_{\al}\eta_{\al}\ox
\et_{\al}\ox v).
$$
The traceless piece of the sum
is the traceless part of  the first summand, while the trace part
of  the sum is
$$
(\frac{2}{m}(w+1)-1)\langle Z,Y\rangle (\sum_{\al}\eta_{\al}\ox\et_{\al}\ox
v).
$$
Consequently, $\pi_1\circ(\na^{\om})^2$ is an invariant operator
for $w=-1,$ $\pi_2\circ(\na^{\om})^2$ is invariant for
$w=\frac{m-2}{2}$ and
$\pi_3\circ(\na^{\om})^2$ is invariant for any value of $w$.

We can now compute the form of those three invariant operators.

{\bf 2.}
Let $\la'=0, w=-1,$ let $\rho'=e^1,\si'=2e^1,$ so
$\pi^{\la\rh\si}=\pi_1.$
Note that
$[e_k,e^i]=e_k^i-e_i^k+\de_k^i\Bbb I_m;$
the semisimple part of $\goth g_0$ is acting
trivially and
$$
\la([X,\Ga .Y])s=(-w)\langle\Ga Y,X\rangle s.
$$

Hence the invariant operator can be written as
$$
\pi_1[(\nabla^{\ga}_a\nabla^{\ga}_b+\Ga_{ab})s]=
[\nabla^{\ga}_{(a}\nabla^{\ga}_{b)_0}+\Ga_{(ab)_0}]s,
$$
where the brackets indicate the symmetrization and the subscript
$0$ means the trace free part.

{\bf 3.}
       Let $\la'=\si'=0,\rh'=e^1; w=\frac{m-2}2$
then
$\la([X,\Ga .Y])s=\frac{2-m}2\sum_{ij}\Ga_{ij}X^jY^i.$
The corresponding projection $\pi_2$ is here just the trace and
 we can express the operator $\pi_2\circ(\na^{\om})^2$
 in a more standard form using
$$
\operatorname{Tr}\left\{
\frac{2-m}2\left[\frac{-1}{m-2}\left(
R_{ab}-\frac{\de_{ab}}{2(m-1)}R\right)\right]\right\}=
\frac{m-2}{4(m-1)}R,
$$
where we used formula \nmb!{6.3}.(1).
Hence we get  the conformally invariant Laplace operator
$$
\pi_2\circ(\na^{\om})^2=
g^{ab}\nabla^{\ga}_a\nabla^{\ga}_b+\frac{m-2}{4(m-1)}R.
$$
This is an example of a so called nonstandard operator.

{\bf 4.}
Let $\la'=0,\rh'=e^1,\si'=e^1+e^2.$ Then $\pi^{\la\rh\si}$ is the projection
to $\La^2(\gog_{-1}^*)\ox V_0;$ i.e. the antisymmetrization.
The tensor $\Ga $ is symmetric, so
$$
\pi_3\circ(\na^{\om})^2=
\na^{\ga}_{[a}\na^{\ga}_{b]}.
$$
Hence we have got a zero order operator  in this case given
by the action of the curvature. In  this case, however,
it is the trivial operator, due to the fact that
the action of $\gog_0$ on $V_0$ is trivial.
But it shows a possibility that for more complicated
representations $V_{\la},$ (e.g. for one forms), we could
get in such a way nontrivial zero order action by the curvature.

{\bf 5.}
To have a more complicated example, let us consider a simple
third order operator. Take $\la'=0,\rh'=e^1,\si'=2e^1$ and
$\ta'=3e^1.$           The projection $\pi^{\la\rh\si\ta}$ is uniquely
defined
by iterated projections to factors having the corresponding
highest weights. The projection $\pi^{\la\rh\si\ta}$ is
the projection to the traceless part of the third symmetric
power.

We can now repeate the computation described in Example 1.
The projection $\pi^{\la\rh\si\ta}$ factorizes through
the projection to $S^3(\gog_{-1}^*)\ox V_{\la},$ hence
the order of the factors is irrelevant, moreover $\pi^{\la\rh\si\ta}$
kills all trace terms. Hence all elements used in Lemma 5.9 have
the form
$$
3(w+2)Z\ox Y_1\ox Y_2\ox v.
$$
The choice $w=-2$ anihilates them all and for this value of
conformal weight, we get a conformally invariant operator.

Proposition \nmb!{4.8} describes the form of the correction terms. Due to the
fact that action of the orthogonal group is trivial, we get for
the first term
$$
\la([X,(\nabla^{\ga}_Z\Ga).Y])s
=2(\nabla^{\ga}_{(a}\Ga^{\strut}_{bc)_0})s.
$$
The next two terms
$\la([X,\Ga.Y])(\nabla^{\ga}_Zs)+
\la([X,\Ga.Z])(\nabla^{\ga}_Ys)$
 lead (due to symmetrization) to the term
$4\Ga^{\strut}_{(ab}\nabla^{\ga}_{c)_0}s$.
The last term
$$
(\la^{(1)}([Y,\Ga.Z])((\nabla^\ga s)(X)
$$
can be written as
$
\la([Y,\Ga.Z])(\nabla^\ga s)(X)+\nabla^\ga{[X,[Y,\Ga.Z]]}s.
$

Using
$[X,[Y,\Ga.Z]]=\langle X,Y\rangle\Ga.Z-\langle X,\Ga.Z\rangle Y-\langle\Ga.Z,Y\rangle X,$
we see that the first term on the right hand side will disappear due
to
the projection to the traceless part and the other two will
cancel the contribution coming from the previous term $\la(...).$
Hence we get the operator
$$
\nabla^\ga_{(a}\nabla^\ga_b\nabla^\ga_{c)_0}s + 4
\Ga^{\strut}_{(ab}\nabla^\ga_{c)_0}s + 2(\nabla^\ga_{(a}\Ga^{\strut}_{bc)_0})s.
$$

The examples shown above illustrate  possibilities of
our approach to construct and to compute the form of
invariant operators. To make computation effective for
a general representations, it is necessary to
use appropriate Casimir operators.
In the next part of the series, we shall use this approach
to describe explicitely the broad family of the so called
standard operators for all AHS structures.

\enddocument